\documentclass{amsart}
\usepackage[utf8]{inputenc}
\usepackage{mathtools, amsmath, amssymb, amsthm}
\usepackage{mathabx}

\usepackage{hyperref}
\usepackage{cleveref}

\usepackage{csquotes}
\usepackage{tikz-cd}

\usepackage{tikz}
\usetikzlibrary{arrows, automata, positioning, patterns, snakes}

\usepackage{circuitikz} 
 
\newtheorem{theorem}{Theorem}[section]
\newtheorem{corollary}[theorem]{Corollary}
\newtheorem{lemma}[theorem]{Lemma}

\newtheorem{proposition}[theorem]{Proposition}
\newtheorem{conjecture}[theorem]{Conjecture}

\newtheorem{theoremA}{Theorem}[section]

\newtheorem{corollaryA}[theoremA]{Corollary}

\newtheorem{propositionA}[theoremA]{Proposition}

\newtheorem{customquestionA}{Question}

\theoremstyle{definition}
\newtheorem{definition}[theorem]{Definition}

\newtheorem*{lemma*}{Lemma}
\newtheorem*{proposition*}{Proposition}
\newtheorem*{theorem*}{Theorem}
\newtheorem*{corollary*}{Corollary}

\theoremstyle{remark}
\newtheorem{remark}[theorem]{Remark}

\theoremstyle{remark}

\theoremstyle{definition}

\newcommand{\sCop}{\operatorname{sCop}}
\newcommand{\wCop}{\operatorname{wCop}}
\newcommand{\Cay}{\operatorname{Cay}}

\newcommand{\rhoZ}{(4\rho\mathbb{Z})^2}
\newcommand{\Zo}{\mathbb{Z}_o^2}
\newcommand{\Z}{\mathbb{Z}} 
\newcommand{\N}{\mathbb{N}}

\newcommand{\into}{\hookrightarrow}

\newcommand{\floor}[1]{\lfloor #1 \rfloor}

\renewcommand{\leq}{\leqslant}
\renewcommand{\geq}{\geqslant}

\begin{document}
\title[Cops and robbers for hyperbolic and virtually free groups]{Cops and robbers for hyperbolic\\ and virtually free groups}

\author{Raphael Appenzeller and Kevin Klinge}

\date{\today}      
\address{Department of Mathematics, Universität Heidelberg, Germany}
\email{rappenzeller@mathi.uni-heidelberg.de}

\address{Department of Mathematics, Universität Heidelberg, Germany}
\email{kklinge@mathi.uni-heidelberg.de}

\def\subjclassname{\textup{2020} Mathematics Subject Classification}
\expandafter\let\csname subjclassname@1991\endcsname=\subjclassname

\begin{abstract}
Lee, Martínez-Pedroza and Rodríguez-Quinche define two new group invariants, the strong cop number $\sCop$ and the weak cop number $\wCop$, by examining winning strategies for certain combinatorial games played on Cayley graphs of finitely generated groups. We show that a finitely generated group $G$ is Gromov-hyperbolic if and only if $\sCop(G) = 1$. We show that $G$ is virtually free if and only if $\wCop(G)=1$, answering a question by Cornect and Martínez-Pedroza. We show that $\sCop(\mathbb{Z}^2) = \infty$, answering a question from the original paper. It is still unknown whether there exist finite cop numbers not equal to 1, but we show that this is not possible for CAT(0)-groups. We provide machinery to explicitly compute strong cop numbers and give examples by applying it to certain lamplighter groups, the solvable Baumslag-Solitar groups, Grigorchuk's group and Thompson's group $F$.

\end{abstract}

\maketitle

\section{Introduction}

Cops and Robbers is a class of combinatorial two-player games that are usually played on graphs or more generally, on metric spaces. Lee, Mart\'inez-Pedroza and Rodr\'iguez-Quinche \cite{LMR23} introduced a version where the existence of a winning strategy for either player is invariant under quasi-isometries, making the game interesting from a geometric group theory perspective. The game depends on a parameter called the \emph{number of cops}. The cop number of a graph is the smallest number of cops such that there exists a strategy for the cop player where they always win. The cop number of a finitely generated group is the cop number of one of its Cayley graphs.
In this work, we characterize the cop numbers in terms of well-known geometric invariants and resolve some open questions in the literature.

There are two versions of the game: the strong cop game and the weak cop game, yielding two notions of cop numbers, which are known to be different in general. Hyperbolic groups are known to have strong cop number \(1\) \cite[Theorem A]{LMR23}. We show that the converse also holds. Thus the strong cop number provides a new characterization of hyperbolic groups.
\begin{theoremA}[\Cref{thm:scop-hyperbolic}]\label{thmA:scop-hyperbolic}
	A finitely generated group \(G\) is Gromov-hyperbolic if and only if its strong cop number is \(1\).
\end{theoremA}

We provide a similar characterization for the weak cop number. As the strong cop number is at most equal to the weak cop number, we know that the weak cop number can be $1$ only for hyperbolic groups.
However, \cite[Corollary J]{LMR23} states that one-ended hyperbolic groups have infinite weak cop number.
We prove the following theorem that was first conjectured in \cite[Question 1.4]{CoMa24arxiv}.
\begin{theoremA}[\Cref{thm:wcop-virtfree}]\label{thmA:wcop-virtfree}
	A finitely generated group $G$ is virtually free if and only if its weak cop number is $1$.
\end{theoremA}

The weak cop numbers of some non-hyperbolic groups have been calculated. We upgrade some of these results to the strong cop numbers. The weak cop number of $\mathbb{Z}^2$ is $\infty$ \cite[Theorem C]{LMR23}, but the value of the strong cop number of $\mathbb{Z}^2$ was left as an open question \cite[Question D]{LMR23}, which we resolve.
\begin{theoremA}[\Cref{thm:scop-Zn}]\label{thmA:scop-Zn}
	For any \(n \geq 2\), \(\Z^n\) has infinite strong cop number.
\end{theoremA}
By \cite[Theorem E]{LMR23}, groups that retract to $\mathbb{Z}^2$, such as Thompson's group $F$, therefore also have strong cop number $\infty$. Another consequence of \Cref{thmA:scop-Zn} is the following.
\begin{theoremA}[Therems \ref{thm:scop-products} and \ref{thm:scop-Grigorchuk}]\label{thmA:scop-Grigorchuk}
    For finitely generated infinite groups $G$, $G\times G$ has infinite strong cop number. Grigorchuk's group and all other finitely generated branch groups have infinite strong cop number.
\end{theoremA}
Cornect and Mart\'inez-Pedroza show that the weak cop number of lamplighter groups is $\infty$ \cite[Theorem 1.5]{CoMa24arxiv}. When the lamp group $L$ is finite, we strengthen this result for the strong cop numbers.
\begin{theoremA}[\Cref{thm:scop-lamplighter}]\label{thmA:scop-lamplighter}
    Let \(L\) be a non-trivial finite group. Then the lamplighter group $L \wr \mathbb{Z}$ has infinite strong cop number.
\end{theoremA}
	
Furthermore, we find the cop numbers of the solvable Baumslag-Solitar groups.
\begin{theoremA}[\Cref{thm:scop-bs}]\label{thmA:scop-bs}
    The solvable Baumslag-Solitar groups $\operatorname{BS}(1,n)$ have infinite strong cop number.
\end{theoremA}
An important concept in the proof of Theorems \ref{thmA:scop-hyperbolic}, \ref{thmA:scop-Zn}, \ref{thmA:scop-lamplighter} and \ref{thmA:scop-bs} is what we call a \emph{meta-stage}:
where a \emph{stage} consists of the cop player taking a move followed by the robber player taking a move, a meta-stage consists of multiple moves for both players.
These moves have to be taken alternately as per the rules of the game.
But by choosing them carefully, we can pretend that first one player takes all their moves within a single meta-stage and then the other player gets their moves, without providing a game-changing advantage to either player. From this point of view, results about the weak cop game can often be upgraded to the strong cop game. The existence of a homothety is one condition that enables such an upgrade. We call this strategy meta-gaming.

\begin{theoremA}[Meta-gaming theorem, for details see \Cref{thm:meta-gaming} and \Cref{cor:self-meta-gaming}]\label{thmA:meta-gaming}
If there is a sequence of quasi-surjective quasi-homotheties from a graph $\Delta$ to a graph $\Gamma$, then the strong cop number of $\Gamma$ is equal to the weak cop number of $\Delta$.
\end{theoremA}

\begin{figure}[ht]
\centering
\includegraphics[width=0.8\textwidth]{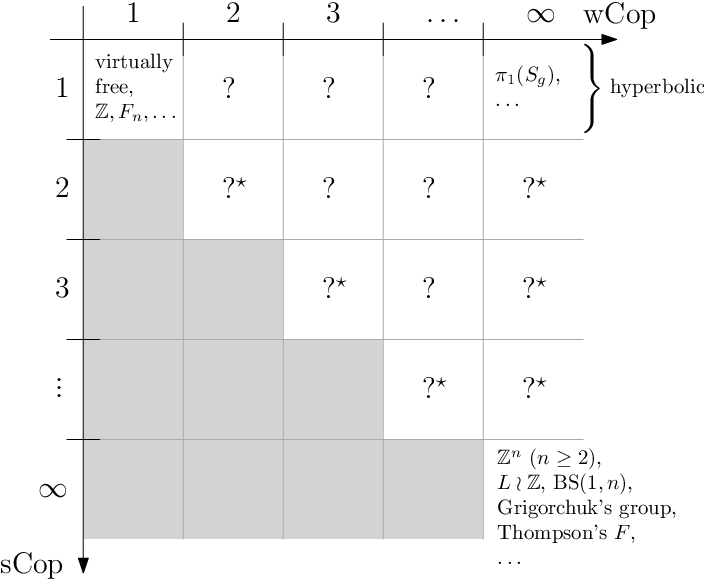} 
\caption{An overview of known results for the strong and weak cop numbers. The existence of groups with combinations of cop numbers marked with ?\ is an open question. The marking $?^\star$ indicates that no group is known, but graphs are known.}
\label{fig:overview}
\end{figure}

In fact, Theorems \ref{thmA:scop-Zn} and \ref{thmA:scop-lamplighter} are obtained from the meta-gaming theorem. On the other hand, when the weak and the strong cop numbers of a group are not equal, the meta-gaming theorem implies that homotheties cannot exist. By Theorems \ref{thmA:scop-hyperbolic} and \ref{thmA:wcop-virtfree}, examples of such groups are hyperbolic groups that are not virtually free. The following are consequences for such groups and spaces that do not involve cop numbers. While the statements may be well known to experts, they serve as a proof of concept that the cop numbers can have broader applications.
\begin{corollaryA}[\Cref{cor:no-homothety-in-hyp-groups}]\label{corA:no-homothety-in-hyp-groups}
If $\Gamma$ is the Cayley graph of a finitely generated hyperbolic group that is not virtually free, then there does not exist a sequence of quasi-surjective homotheties.
\end{corollaryA}
\begin{corollaryA}[\Cref{cor:no-homothety-in-hyp-plane}]\label{corA:no-homothety-in-hyp-plane}
     The hyperbolic plane $\mathbb{H}^2$ does not admit a surjective homothety $h \colon \mathbb{H}^2 \to \mathbb{H}^2$.
\end{corollaryA}

The main remaining open question is whether there exist groups with intermediate cop numbers.
\begin{customquestionA}[\cite{LMR23}, Question K]\label{qA:intermediate}
Are there any finitely generated groups with weak and/or strong cop numbers not equal to $1$ or $\infty$?
\end{customquestionA}

If the answer to Question \ref{qA:intermediate} is positive, then the cop numbers could become a valuable tool in geometric group theory since they would serve as close generalizations of being hyperbolic and/or being virtually free. If the answer is negative, Theorems \ref{thmA:scop-hyperbolic} and \ref{thmA:wcop-virtfree} fully classify the cop numbers in terms of known invariants. In the special case of CAT(0)-groups we show that there are no intermediate strong cop numbers.
\begin{propositionA}[Proposition \ref{prop:intermediate-CAT0}]\label{propA:intermediate-CAT0}
    If $G$ is a CAT(0)-group, then $\sCop(G) \in \{ 1,\infty \} $.
\end{propositionA}
Figure \ref{fig:overview} gives an overview of the results proven here.

\medskip
\textbf{Outline.}
We start by defining the cop numbers and recalling some results about them in \Cref{sec:def}. In \Cref{sec:hyp} we prove the characterizations of hyperbolic and virtually free groups in terms of their cop numbers, Theorems \ref{thmA:scop-hyperbolic} and \ref{thmA:wcop-virtfree}. Section \ref{sec:meta} is dedicated to metagaming with homotheties. We start by showing that $\sCop(\mathbb{Z}^2) = \infty$ (\Cref{thmA:scop-Zn}) before generalizing the argument to the meta-gaming theorem (\Cref{thmA:meta-gaming}). In Section \ref{sec:applications}, we collect the results about non-existence of homotheties (Corollaries \ref{cor:no-homothety-in-hyp-groups}, \ref{cor:no-homothety-in-hyp-plane}), Lamplighter groups (\Cref{thmA:scop-lamplighter}), solvable Baumslag-Solitar groups (\Cref{thmA:scop-bs}), Grigorchuk's group (Theorem \ref{thmA:scop-Grigorchuk}) and Thompson's group $F$. In Section \ref{sec:intermediate} we argue that CAT(0)-groups do not admit intermediate cop numbers.

\medskip
\textbf{Acknowledgements.}  

We would like to thank Francesco Fournier-Facio and Luca de Rosa, with whom we played many rounds of Cops and Robbers and had many fruitful discussions. For discussions about hyperbolicity and related questions, we would also like to thank Tommaso Goldhirsch, Hjalti Isleifson, Alexander Lytchak, Bianca Marchionna, Jos\'e Pedro Quintanilha and Claudio Llosa Isenrich. We are thankful for a detailed report by an anonymous referee.
The authors were partially supported by RTG 2229 ``Asymptotic Invariants and Limits of Groups and Spaces'' funded by Deutsche Forschungsgemeinschaft (DFG, German Research Foundation).
During the publication process of this paper, some of our results were independently obtained by \cite{EGG25arxiv} using different methods. There has also been partial progress on Question~\ref{qA:intermediate}: \cite{Leh25arxiv} shows that there are no intermediate weak cop numbers.

\section{The weak and the strong cop numbers}\label{sec:def}

We start by recalling the definition of the cops and robbers game by Lee, Martínez-Pedroza, Rodríguez-Quinche \cite{LMR23}.
The two-player game with perfect information is defined on a connected undirected graph $\Gamma$ without double edges and loops.

Apart from the graph itself, the game depends on six parameters:
\begin{enumerate}
	\item the \emph{number of cops} \(n\in \N \coloneq \Z_{\geq 1}\),
	\item \(v\) the \emph{treasure} or \emph{base point}: a vertex in \(\Gamma\),
	\item \(\sigma \in \N\) the \emph{cop speed},
	\item \(\psi \in \N\) the \emph{robber speed},
	\item \(\rho \in \N\) the \emph{cop reach} and
	\item \(R \in \N\) the \emph{robber radius}.
\end{enumerate}
The game is played as follows. First, the cop player chooses $n$ vertices of $\Gamma$ as starting positions.
Then the robber chooses one vertex of $\Gamma$ as their starting position. The players take turns moving, starting with the cops.
In each turn, the cop-player can move each cop along a path of length at most $\sigma$ to some vertex in $\Gamma$.
The robber-player then moves the robber along a path of length at most $\psi$ ending in a vertex. One cop-move followed by one robber-move is called a \emph{stage}.

The \emph{cop-player wins} if the \emph{robber is captured}: during any stage in the game the robber's position is at distance at most $\rho$ from some cop's position; or if the cops can \emph{eventually defend $v$}: after some finite initial time, the robber never enters the closed ball of radius $R$ around $v$.
Note that the robber not only needs to be at least $\rho$ away from the cops at the beginning and ends of stages but also during the whole movement of the robber. The robber needs to be able to take their path edge by edge without ever coming $\rho$-close to a cop.

The \emph{weak cop number} $\wCop(\Gamma)$ of $\Gamma$ is the smallest value $n$ such that for any $v\in\Gamma$, the cop-player wins the game when first the cop-player chooses their speed $\sigma$ and reach $\rho$ and then the robber-player chooses their speed $\psi$ and the radius $R$. Formally
\begin{align*}
    \wCop(\Gamma) = \inf \left\{  n\in \mathbb{N} \colon  \begin{array}{l}
        \forall v \in \Gamma, \exists \sigma \in \mathbb{N}, \exists \rho \in \mathbb{N}, \forall \psi \in \mathbb{N}, \forall R \in \mathbb{N} \colon \\
        \exists\text{ a winning strategy for the cop-player with $n$ cops.}
    \end{array}
    \right\}.
\end{align*}
The \emph{strong cop number} $\sCop(\Gamma)$ of $\Gamma$ is defined similarly, except that the cop-player chooses $\rho$ after the robber player chooses $\psi$, formally
\begin{align*}
    \sCop(\Gamma) = \inf \left\{  n\in \mathbb{N} \colon  \begin{array}{l}
        \forall v \in \Gamma, \exists \sigma \in \mathbb{N}, \forall \psi \in \mathbb{N},  \exists \rho \in \mathbb{N}, \forall R \in \mathbb{N} \colon \\
        \exists\text{ a winning strategy for the cop-player with $n$ cops.}
    \end{array}
    \right\}.
\end{align*}
Since $\wCop$ and $\sCop$ are preserved under quasi-isometry \cite[Corollary G]{LMR23}, the cop numbers $\wCop(\Gamma)$ and $\sCop(\Gamma)$ of a Cayley graph $\Gamma = \Cay(G, S)$ of a finitely generated group $G$ do not depend on the finite generating set \(S\).
This allows us to define \(\wCop(G) \coloneq \wCop(\Gamma)\) and \(\sCop(G) \coloneq \sCop(\Gamma)\).
It is an exercise to convince oneself that $\sCop(\Gamma) \leq \wCop(\Gamma)$. Moreover $\wCop(\mathbb{Z}) = \sCop(\mathbb{Z}) = 1$ and this generalizes to infinite trees, such as Cayley graphs of free groups.

\begin{remark}
	The existence of a winning strategy is always independent of the vertex \(v\):
	Suppose that for \(v \in \Gamma\) there exists a winning strategy for the robber using speed \(\rho\) and reach \(R\).
	Then there also exists a winning robber-strategy for \(v' \in \Gamma\) by using the same speed, doing the same moves and picking the reach \(R' \coloneq R + d_\Gamma(v, v')\).
	In the case where \(\Gamma\) is a Cayley graph, we will for simplicity always choose \(v = 1\).
\end{remark}
\begin{remark}
	As the cop numbers are invariant under quasi-isometry and every geodesic metric space is quasi-isometric to a graph, we may define the cop-number of arbitrary geodesic metric spaces.
	It is also possible to define the cops and robbers game directly on geodesic metric spaces where the players move along paths with lengths bounded by their speeds and thus define the cop number of a geodesic metric space. As in \cite[Theorem 4.2]{LMR23} the cop numbers for geodesic metric spaces are invariant under quasi-isometries, so these two methods yield the same definition.
\end{remark}

The following theorem is an important result about retractions, that was used by \cite{LMR23} to establish that the cop numbers are invariant under quasi-isometry.
\begin{theorem}[\cite{LMR23} Theorem E]\label{thm:retract}
    Let $\Gamma, \Delta$ be graphs. If there is a quasi-retraction from $\Gamma$ to $\Delta$, then
    $$
    \wCop (\Delta) \leq \wCop(\Gamma) \quad \text{and} \quad \sCop (\Delta) \leq \sCop(\Gamma).
    $$
\end{theorem}

\section{Cop numbers and known geometric invariants.}\label{sec:hyp}
\subsection{Strong cop number and hyperbolicity}

One of the most important invariants in geometric group theory is hyperbolicity. A geodesic metric space is \emph{hyperbolic} if there exists $\delta$ such that all triangles are $\delta$-slim (if the triangle is given by geodesics $[x,y], [y,z], [x,z]$, then $[x,z] \subseteq N_\delta([x,y] \cup [y,z]) $). A finitely generated group is \emph{hyperbolic} if its Cayley graphs are. We use the following alternative characterization of hyperbolicity for graphs in terms of thin bigons due to \cite{Pap95}. A \emph{bigon} in a graph $\Gamma$ is a pair of geodesic paths $\gamma, \gamma' \colon [0, \ell]\cap \mathbb{Z} \to \Gamma$ with $\gamma(0) = \gamma'(0)$ and $\gamma(\ell) = \gamma'(\ell)$. For $\delta \in \mathbb{R}$, we say that bigons in a graph $\Gamma$ are \emph{$\delta$-thin} if for all $t \in [0,\ell] \cap \mathbb{Z}$, $d(\gamma(t),\gamma'(t) ) \leq \delta$.


\begin{theorem}\cite[Theorem 1.4]{Pap95}\label{thm:papasoglu}
    Let $G$ be a finitely generated group and $\Gamma$ the Cayley graph associated to a finite generating set of $G$. If there is a $\delta$ such that all bigons are $\delta$-thin, then $G$ is hyperbolic.
\end{theorem}

\begin{remark}
The theorem is only stated for Cayley graphs of groups in \cite{Pap95}, but the proof given there is valid for graphs in general. One should however be careful and note that the theorem does not hold for general geodesic metric spaces, such as $\mathbb{R}^2$. For general geodesic metric spaces, Pomroy proved a characterization of Gromov-hyperbolicity in terms of quasi-geodesic bigons \cite[Theorem 22]{ChNi07}, but for us, it is easier to work with the characterization given here.
\end{remark}
The converse of Theorem \ref{thm:papasoglu} is also true: in a hyperbolic graph, every bigon is a degenerate triangle and $\delta$-slim bigons ($\gamma \subseteq N_\delta(\gamma')$) are $2\delta$-thin, see for instance \cite[Lemma 4]{ChNi07}.
\begin{corollary}\label{cor:papasoglu}
    A graph is hyperbolic if and only if there is a $\delta$ such that all bigons are $\delta$-thin.
\end{corollary}
Lee, Martínez-Pedroza and Rodríguez-Quinche showed that every hyperbolic graph $\Gamma$ has strong cop number $1$ \cite[Theorem A]{LMR23}. The following theorem shows that this is a characterization of hyperbolicity.
\begin{theorem}\label{thm:scop-hyperbolic}
    A graph $\Gamma$ is hyperbolic if and only if $\sCop(\Gamma) = 1$.
\end{theorem}
\begin{proof}
    Because of \cite[Theorem A]{LMR23}, it suffices to show that if $\Gamma$ is not hyperbolic, then $\sCop(\Gamma)>1$. Let $\Gamma$ be non-hyperbolic, then by Corollary \ref{cor:papasoglu}, for every $\delta >0 $ there exist bigons that are not $\delta$-thin. We now describe a winning strategy for the robber against one strong cop. The game starts with the cop player choosing their speed $\sigma$. The robber then chooses speed
    $$
    \psi = 96 \sigma.
    $$
		The cop then chooses reach $\rho$, which we may assume without loss of generality to satisfy $\rho \geq \sigma$. Since $\Gamma$ is not hyperbolic, there exists a $\delta > 32\rho$ and a bigon $(\gamma, \gamma')$ that is $\delta$-thin, but not $(\delta-1)$-thin. The robber now chooses a radius $R$ large enough so that a $\delta$-neighborhood of the bigon is contained in the ball of radius $R$ around the base point. Next, the cop chooses a starting position. The distance $\lambda := \delta/16 > 2 \rho \geq \rho + \sigma$ will play an important role in the following, it is a distance at which the robber is safe from the cop, even after the cop moves once. The bigon $(\gamma, \gamma')$ has the property that there exists $t \in [0,\ell]$ such that $d(\gamma(t),\gamma'(t)) = \delta = 16\lambda$ and we denote $p := \gamma(t)$ and $p':=\gamma'(t)$.
    At least one of the two points $p,p'$ has distance larger than $5 \lambda$ from the starting position of the cop, and the robber chooses that point as their starting point $r_0$. Without loss of generality $r_0 = p$. The game now starts with a move from the cop.

    The strategy of the robber during the game is as follows. As long as the distance to the cop is at least $5\lambda$, the robber just stays where they are and waits. As soon as the cop gets closer than $5\lambda$ to the robber, the robber will move from $p$ to $p'$ along one of the two paths $\eta_{\pm}$ starting at $p$, following along $\gamma$ to $\gamma(t \pm \delta)$, then along some geodesic to $\gamma'(t \pm \delta)$ and finally along $\gamma'$ to $p'$; see \Cref{fig:scophyp_main} for an illustration. In case $t+\delta > \ell$ or $t-\delta < 0$, we let $\eta_{\pm}$ just consist of the paths via $\gamma(0)$ or $\gamma(\ell)$ without the middle shortcut. Once arrived at $p'$, the robber repeats the strategy of waiting until the cop gets closer than $5\lambda$ and then switches back to $p$. By the choice of $R$, the robber never strays too far from the base point, so to show that this is a winning strategy for the robber, it suffices to show that the robber can choose and follow one of the paths $\eta_+$ or $\eta_-$ without ever being caught.

\begin{figure}[h]
  \centering
  \includegraphics[width=0.8\textwidth]{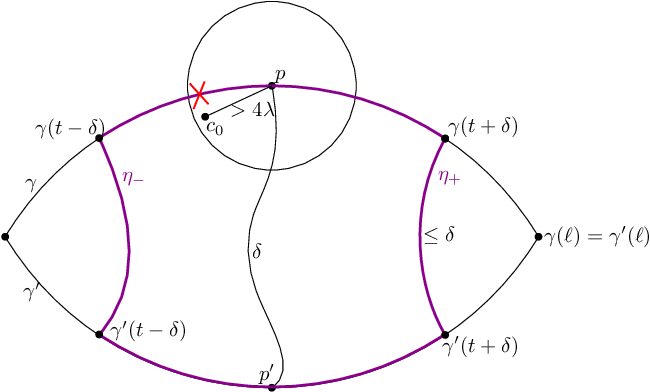} 
  \caption{The robber finds a $\delta$-thin, but not $(\delta-1)$-thin bigon and waits at $p$. Once the cop gets closer than $5\lambda$, the robber moves from $p$ to $p'$ along one of two possible paths $\eta_+$ or $\eta_-$, at most one of which can be blocked by the cop $c_0$. }
  \label{fig:scophyp_main}
\end{figure}
\begin{figure}[h]
  \centering
  \includegraphics[width=1\textwidth]{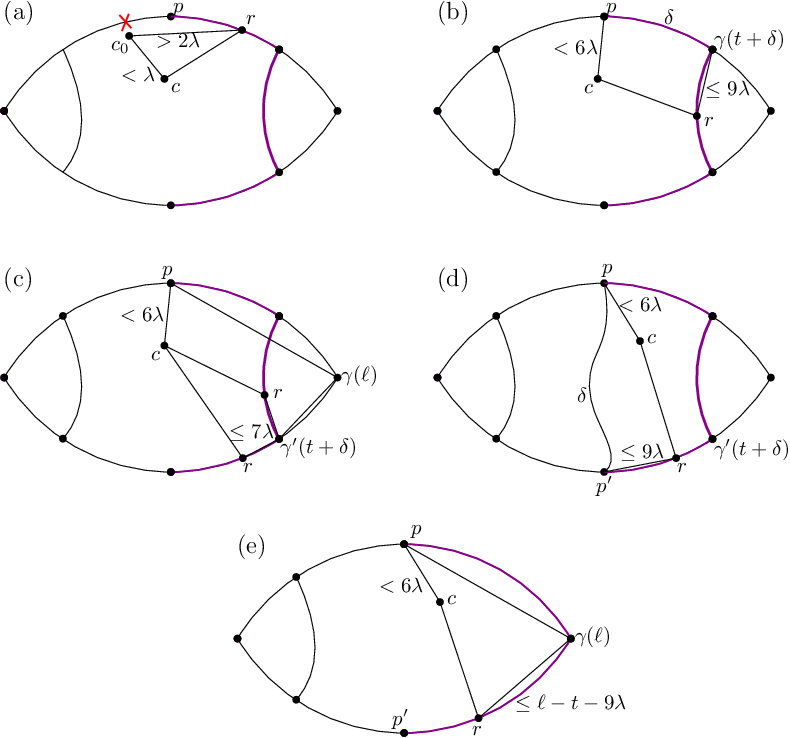} 
  \caption{To prove that the robber $r$ is not caught by the cop $c$ during their movement along $\eta_+$ we consider five cases. Case (e) is only needed if $\ell < t+\delta$.  }
  \label{fig:scophyp_cases}
\end{figure}

    Let $r_0 = p$ be the starting point of the robber and $c_0$ be the position of the cop when they get closer than $5\lambda$ for the first time. Since the cop can move at most $\sigma < \rho < \lambda$, we also have
    \begin{align}
        4 \lambda < d(p,c_0) < 5 \lambda. \label{eq:hyp_pc0}
    \end{align}
    Let $r$ be the position of the robber somewhere on the path $\eta_{\pm}$, and let $c$ be any possible position of the cop during the movement of the robber along $\eta_{\pm}$. Moving along $\eta_{\pm}$ amounts to moving a distance of at most $3\delta$, which takes at most $\lceil 3\delta/\psi \rceil = \lceil \delta/(32 \sigma) \rceil $ steps.
		During this time, the cop can move a distance of at most $ \lceil \delta/(32\sigma) \rceil \sigma \leq \lceil \delta/ 32 \rceil  \lambda/2 $, so $d(c,c_0) \leq \lceil \lambda/2\rceil< \lambda$. Together with the triangle inequality and \Cref{eq:hyp_pc0}, this shows
    \begin{align}
        3 \lambda < d(p,c) < 6 \lambda . \label{eq:hyp_pc}
    \end{align}
		We claim that if for some $s\geq 0$, $d(c_0, \gamma(t-s)) < 2\lambda$, then $s > 2\lambda$ and for every $u >0$, $d(c_0,\gamma(t+u)) > \lambda$. Intuitively the claim states that if $\eta_-$ is blocked, then (at least the part along $\gamma$ of) $\eta_+$ is not blocked. Indeed, the triangle inequality gives $s = d(\gamma(t-s),p)  \geq d(c_0,p) - d(c_0,\gamma(t-s)) > 2\lambda$. If $u< 2\lambda$, then
    $$
    d(c_0,\gamma(t+u)) \geq d(c_0,p) - d(\gamma(t+u),p) > 4\lambda - 2\lambda = 2\lambda.
    $$
    by \Cref{eq:hyp_pc0}. If $u>2\lambda$, then
    $$
    d(c_0,\gamma(t+u)) \geq d(\gamma(t-s),\gamma(t+u)) - d(\gamma(t-s),c_0) > (s+u) - 2\lambda  >2\lambda
    $$
    concludes the proof of the claim.

    Analogously, if for some $s\geq 0$, $d(c_0, \gamma(t+s))<2\lambda$, then for every $u \geq 0$, $d(c_0,\gamma(t-u))> 2\lambda$. This implies that at most one of the paths $\eta_+$ or $\eta_-$ can be blocked. Let us assume without loss of generality that $\eta_+$ is not blocked. For the entire movement of the robber, we consider the following cases, illustrated in Figure
    \ref{fig:scophyp_cases}:

    \begin{enumerate}
        \item[(a)] If $r = \gamma(t + u)$ for some $u\geq 0$, then $d(c,r) \geq d(c_0,\gamma(t+u)) - d(c_0,c) > 2\lambda - \lambda = \lambda$.
        \item[(b)] If $d(r,\gamma(t+\delta)) \leq 9\lambda$, then $d(r,c) \geq d(\gamma(t+\delta), p) - d(\gamma(t+\delta) , r) - d(c,p) > 16\lambda - 9\lambda -6\lambda = \lambda$ by Equation (\ref{eq:hyp_pc}). \label{eq:hyp-2ndcase}
        \item[(c)] If $d(r,\gamma'(t + \delta)) \leq 7\lambda$, then
    \begin{align*}
        d(r,c) &\geq d(\gamma(\ell),p) - d(\gamma(\ell),\gamma'(t+\delta)) - d(r,\gamma'(t + \delta)) - d(c,p) \\
        &> (\ell-t) - (\ell-(t+\delta)) - 7\lambda- 6\lambda = 3 \lambda > \lambda.
    \end{align*} \label{eq:hyp-3rdcase}
    \item[(d)]
    If $d(r,p') \leq 9 \lambda$, then $d(r,c) \geq d(p,p') - d(r,p') - d(p,c) > \delta - 9\lambda - 6 \lambda = \lambda$.
    \end{enumerate}
    If $t+\delta > \ell$, there is no point $\gamma(t+\delta)$. Then we replace cases (b) and (c) by the following.
      \begin{enumerate}
      \item[(e)] If $d(r,p')\geq 9\lambda$ and $r \in \operatorname{im}(\gamma')$, then $d(\gamma(\ell),r) = d(p',\gamma(\ell)) - d(r,p') \leq \ell - t - 9\lambda $ and $d(r,c) \geq d(p,\gamma(\ell)) - d(p,c) - d(\gamma(\ell), r) > (\ell - t) - 6\lambda - (\ell - t - 9\lambda) = 3\lambda > \lambda$.
      \end{enumerate}
		While the robber moves along $\eta_+$, they are always in at least one of the cases above and thus never get caught.  Finally, once the robber reaches $p'$, we have that $d(p',c) \geq d(p',p) - d(p,c) > 10 \lambda > 5\lambda$, so the strategy can repeat and the robber has a winning strategy.
\end{proof}

\subsection{Weak cop numbers and virtually free groups}
In this section we characterize the graphs and groups with weak cop number $1$, answering a question of Cornect and Mart\'inez-Pedroza \cite[Question 1.4]{CoMa24arxiv}.
\begin{theorem}\label{thm:wcop-virtfree}
    A finitely generated group $G$ is virtually free if and only if it satisfies $\wCop(G)=1$.
\end{theorem}
\begin{theorem}\label{thm:wcop-qi-tree}
	Let \(X \) be a geodesic metric space.
	Then \(\wCop(X) = 1\) if and only if \(X\) is a quasi-tree.
\end{theorem}
The two theorems are equivalent since the geometry and algebra of virtually free groups are related by the following well-known theorem.
\begin{theorem}[\cite{DuKa2018} Theorem 20.45, or \cite{KrMo08} Theorem 3.28, or \cite{button23arxiv} Corollary 10.7]\label{thm:qtreevirtfree}
    A finitely generated group is virtually free if and only if its Cayley graph is a quasi-tree.
\end{theorem}
We use the following characterization of quasi-trees.
\begin{theorem}[\cite{Manning2003GeometryOP}, Theorem 4.6, Bottleneck Property]
	A geodesic metric space \(X\) is a quasi-tree if and only if it satisfies the \emph{Bottleneck Property}:

	There is some \(\delta > 0\) such that for all \(x, z \in X\) there exists \(y \in X\) with \(d(x, y) = d(y, z) = \frac12 d(x,z)\) and the property that any path \(\gamma \) from \(x\) to \(z\) satisfies \(d(\gamma , y) \leq \delta\).
\end{theorem}

\begin{proof}[Proof of Theorem \ref{thm:wcop-qi-tree}]

    \begin{figure}
  \centering
  \includegraphics[width=0.5\textwidth]{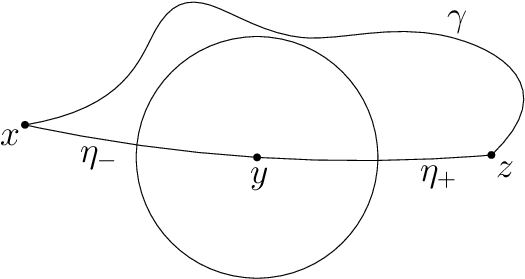} 
  \caption{ The negation of the Bottleneck property guarantees the existence of a path $\gamma$ connecting $x$ and $y$ outside a ball of size $6\lambda$ around $y$. }
  \label{fig:wcop-quasitree}
\end{figure}

	If \(X\) is a quasi-tree, then it has \(\wCop(X) = 1\). So suppose that \(X\) is not a quasi-tree.
	We provide a winning strategy for the robber against one weak cop.

	The game starts with the cop player picking speed \(\sigma\) and reach \(\rho\).
	Set \(\lambda \coloneq \sigma + \rho\).
	As \(X\) is not a quasi-tree, it does not satisfy the bottleneck property.
	Hence in \(X\) there are points \(x, y, z\) and a path \(\gamma \) connecting \(x\) to \(z\) such that \(d(x, y) = d(y, z) = \frac12 d(x, z)\) and \(d(y, \gamma) > 6\lambda\).

	Let \(\eta\) be a geodesic from \(x\) to \(z\) through \(y\) and denote by \(\eta_-\) and \(\eta_+\) the parts of \(\eta\) from \(x\) to \(y\) and from \(y\) to \(z\) respectively, see Figure \ref{fig:wcop-quasitree}.
	The strategy for the robber is to end every move in either \(x\) or \(y\) and to move between the two points using either \(\eta_-\) or \(\eta_+\) and \(\gamma\).
	Thus we pick the robber's speed \(\psi\) and reach \(R\) large enough to allow these moves.
	Our strategy is to maintain the following invariant after every move of the robber:
	Let \(c\) denote the cop's position at that time.
	\begin{enumerate}
		\item Either the robber is positioned at \(y\) and \(d(c, y) \geq 4\lambda\)
		\item or the robber is positioned at \(x\) and \(d(c, y) < 4\lambda\).
	\end{enumerate}
	Note that as \[
		d(x, y) \geq d(\gamma , y) > 6\lambda,
	\] this ensures that the robber cannot be caught during any of the cop player's turns. Hence it is enough to prove that this invariant can indeed be maintained.
	Further, there is a point among \(x\) and \(y\) which may serve as a starting position for the robber to satisfy the invariant.

	To describe the robber's moves, suppose first that the robber is currently positioned at \(y\) and hence the cop is at least \(4\lambda\) away from \(y\).
	As long as the cop does not come within distance \(4\lambda\) of \(y\), the robber stays at \(y\).
	If the cop ever moves to some point \(c\) such that \(d(c, y) < 4\lambda\), the robber either moves along \(\eta_-\) to \(x\), or along \(\eta_+\) to \(z\) and then along \(\gamma \) to \(x\). Note that moving along either path takes the robber only one turn.
	Hence it suffices to show that the cop cannot block both paths at the same time.
	That is, we show that $d(\gamma, c) \geq \lambda > \rho$ and \[
		 d(\eta_+, c) \geq \lambda \quad\text{or}\quad d(\eta_-, c) \geq \lambda.
	\]
	By the triangle inequality, \[
		d(\gamma, c) \geq d(\gamma , y) - d(y, c) > 6\lambda - 4\lambda = 2\lambda,
	\]
	so we only need to show that one of \(d(\eta_\pm, c) \geq \lambda\).
	Suppose that \(d(\eta_+, c) < \lambda\) and let \(c' \in \eta_+\) such that \(d(c', c) < \lambda\).
	As the cop just moved into the ball of radius \(4\lambda\) around \(y\), they moved at most \(\sigma \leq \lambda\) into that ball. That is, \(3\lambda \leq d(c, y) < 4\lambda\).
	We want to show that \(d(\eta_-, c) \geq \lambda\).
	As \(\eta\) is a geodesic containing \(c'\), we have that \(d(\eta_-, c') = d(y, c')\) and hence \begin{align*}
		d(\eta_-, c) &~\geq d(\eta_-, c') - d(c, c') = d(y, c') - d(c, c')\\
		&~\geq d(y, c) - d(c, c') - d(c, c')
		\geq 3\lambda - \lambda - \lambda = \lambda> \rho,
	\end{align*}
    finishing the proof that at least one of \(\eta_\pm\) is not blocked by the cop.

	Now suppose that the robber is currently positioned at \(x\).
	As long as the cop stays within the ball of radius \(4\lambda\) around \(y\), the robber stays at \(x\).
	As soon as the cop moves to some point \(c\) such that \(d(c, y) \geq 4\lambda\), we show that the path \(\gamma\), as well as at least one of \(\eta_\pm\) is not blocked by the cop, so the robber can move to \(y\).
	The proof uses essentially the same idea as the first part.

	Again, as the cop just moved to \(c\) from some point that had distance at most \(4\lambda\) from \(y\) and they move only distance \(\sigma \leq \lambda\), we know that \(4\lambda \leq d(c, y) \leq 5\lambda\) and hence \[
		d(\gamma, c) \geq d(\gamma, y) - d(y, c) \geq 6\lambda - 5\lambda = \lambda> \rho.
	\]

	Suppose that \(\eta_+\) is blocked by the cop. That is, there is some \(c' \in \eta_+\) such that \(d(c, c') \leq \sigma < \lambda\).
	Then \[
		d(\eta_-, c) \geq d(\eta_-, c') - d(c', c) = d(y, c') - d(c', c) \geq d(y, c) - 2d(c', c) \geq 4\lambda - 2\lambda > \lambda,
	\] finishing the proof.
\end{proof}

\section{Metagaming with quasi-homotheties}\label{sec:meta}
\subsection{Strong cop number of $\mathbb{Z}^n$.}

The square grid $\Gamma = \mathbb{Z}^2$ has $\wCop(\mathbb{Z}^2) = \infty$ \cite[Theorem C]{LMR23}, but the value of the strong cop number is left as an open question \cite[Question D]{LMR23}. We show that $\sCop(\mathbb{Z}^2)=\infty$ by reducing to the weak cop game.

\begin{theorem}\label{thm:scop-Zn}
 For $n \geq 2$, $\sCop (\mathbb{Z}^n) = \infty$.
\end{theorem}

\begin{proof}
    Since $\mathbb{Z}^n$ retracts to $\mathbb{Z}^2$, it suffices to consider the case $n=2$, by \Cref{thm:retract}.

    We have $\wCop(\mathbb{Z}^2)=\infty$, meaning that there is a winning strategy for the robber against any number of weak cops. In the following, we describe a winning strategy for the robber against any number of strong cops, by reducing to the weak version of the game.

    Let $n\in \mathbb{N}$ arbitrary and let $\sigma_0=2$, $\rho_0=3$. We know that there are $\psi_0, R_0 \in \mathbb{N}$ such that there is a winning strategy for the robber-player. Throughout the proof, the robber-player will access this strategy as an oracle. This copy of $\mathbb{Z}^2$ will be denoted by $\Zo$.

 To show that $\sCop(\mathbb{Z}^2) > n$, we will show that $\exists v \in \mathbb{Z}^2$, $\forall \sigma \in \mathbb{N}, \exists \psi \in \mathbb{N}, \forall \rho \in \mathbb{N}, \exists R \in \mathbb{N}$ such that there is no winning strategy for the cop-player with $n$ cops. We choose $v=0$. Given $\sigma \in \mathbb{N}$, we choose
 $$
 \psi = 4\psi_0 \sigma
 $$
 and given $\rho\in \mathbb{N}$, we choose
 $$
 R = 4 R_0 \rho,
 $$
 where $\psi_0$ and $R_0$ are the values from the weak-cop game described above. We may assume without loss of generality that $\rho$ is a multiple of $\sigma$ by increasing $\rho$.

 The game starts with the cop-player choosing $n$ starting positions $c_1, \ldots, c_n \in \mathbb{Z}^2$. The robber-player chooses their starting position as follows. The robber-player considers the subgraph $\rhoZ \subseteq \mathbb{Z}^2$. Scaling by a factor $4\rho$ gives an isomorphism $\iota \colon \Zo \cong \rhoZ$. For each cop $c_i$ on $\mathbb{Z}^2$, the robber-player considers the location and chooses a closest vertex $p_i$ of $\rhoZ$ (there are at most 4 possible closest vertices). Note that $d(p_i,c_i) \leq 4\rho$. The robber-player now queries the oracle with the cop starting positions $\iota^{-1}(p_1), \ldots, \iota^{-1}(p_n) \in \Zo$ to obtain a starting position in $\Zo$ and via $\iota$ also a starting position $r \in \rhoZ \subseteq \Z$. We note that the robber is not currently captured, as
 \begin{align*}
    12\rho = 4 \rho \cdot \rho_0 < 4\rho \cdot d(\iota^{-1}(p_i),\iota^{-1}(r)) = d(p_i, r) \leq d(p_i,c_i) + d(c_i,r) \leq 4\rho + d(c_i,r)
 \end{align*}
 implies $d(c_i,r) > 8\rho > \rho$.

 During the game, the robber-player continues to play as if the game took place on $\rhoZ$: since the robber's speed $\psi$ is potentially much slower than the grid size $4\rho$ of $\rhoZ$, the robber can not move from one vertex of $\rhoZ$ to another in one stage. Instead, we introduce the notion of a \emph{meta-stage}, which consists of $\rho/\sigma$ many stages on $\mathbb{Z}^2$. In one meta-stage, the robber can move a distance up to
 $$
 \psi \cdot \rho/\sigma = 4 \psi_0 \sigma\cdot \rho/\sigma = 4\rho \psi_0,
 $$
 corresponding to $\psi_0$ many edges on $\Zo$. At the beginning of a meta-stage, when it is the robber-player's turn, they analyze the board and approximate the cop's locations to vertices of $\rhoZ$ as before. The robber-player then queries the oracle, to obtain a path $\gamma$ to a new point in $\Zo$ that prevents the weak cop-player from having a winning strategy. For the rest of the meta-stage, the robber-player ignores the plays of the cop-player and just follows the corresponding path $\iota(\gamma)$ on $\rhoZ \subseteq \Z^2$.
 
We note that during one meta-stage, the cops can move at most
 $$
 \sigma \cdot \rho/\sigma = \rho
 $$
 which is a fourth of an edge length on the graph $\rhoZ$. However, since the cops are not restricted to grid points of $\rhoZ$, their perceived movement on $\rhoZ$ can be larger. Let $c_t,c_{t+1} \in \mathbb{Z}^2$ be a cop's position at the beginning and the end of a meta-stage, and let $p_t,p_{t+1} \in \rhoZ$ be the corresponding pretended locations. Since $d(c_t,p_t) \leq 4\rho$ and $d(c_{t+1},p_{t+1}) \leq 4\rho$,
 $$
 d(p_t, p_{t+1}) \leq d(p_t,c_t) + d(c_t,c_{t+1}) + d(c_{t+1},p_{t+1}) \leq 8\rho+\rho < 3(4\rho)
 $$
 and since distances between graph points on $\rhoZ$ lie in $(4\rho)\mathbb{Z}$, we actually have $d(p_t,p_{t+1})\leq 8\rho$ and thus, from the point of view of $\Zo$, cops can move at most $\sigma_0 = 2$ vertices per meta-stage. The range of the cops is $\rho$ which even together with the $4\rho$ offsets of the actual locations to the perceived locations, is $9\rho \leq 3(4\rho)$, and hence corresponds to at most $\rho_0 = 3$ edges on $\Zo$. Finally, the radius $R = 4R_0\rho$ corresponds to a radius of $R_0$ on $\Zo$.

 From the above discussion, it follows that the movement of the cops and the robber on $\rhoZ$ in terms of meta-stages corresponds to the movement on $\Z^2$ in the weak cop game oracle. From this, it is clear that at the beginning of every meta-stage, the robber is not captured.

 We now argue why the robber never reaches a $\rho$-neighborhood of a cop, even at any point within a meta-stage or a stage. If at any time within a meta-stage, a cop $c \in \Z^2$ is in reach $\rho$ of one of the $4\rho-2$ many points that lie on an edge between two vertices $a,b \in \rhoZ$, where $\iota^{-1}(a)$ is adjacent to $\iota^{-1}(b)$ in $\Zo$, then any pretended location $p$ of $c$ is within distance $5\rho$ of both endpoints $a,b$. Moreover, the position of the cop at the beginning and the end of the meta-stage differs by at most $\rho$. Thus the position of the cop at the beginning and the end of the meta-stage is less than $8\rho$ from both endpoints $a,b \in \rhoZ$. Hence, in the oracle $\Zo$, the corresponding vertices $\iota^{-1}(a),\iota^{-1}(b)$ are reachable by the cop during the whole meta-stage. Hence, in the oracle $\Zo$, the robber will never pass through the vertices $\iota^{-1}(a)$, $\iota^{-1}(b) \in \Zo$, the robber in $\Z^2$ will also never pass through $a,b$ and neither through any points between $a$ and $b$.

 A robber following this strategy can therefore not be captured and the cops cannot eventually defend $v=0$. This means that there is no winning strategy for the cop-player. As this argument works for all $n\in \mathbb{N}$, we have $\sCop{(\mathbb{Z}^2)}= \infty$.
\end{proof}

\subsection{Meta gaming for quasi-homotheties }

The main idea in the proof of Theorem \ref{thm:scop-Zn} is to use a homothety $\mathbb{Z}^2 \to \mathbb{Z}^2$ to play in meta-stages. In this section, we generalize this idea in Theorem \ref{thm:meta-gaming}, where we show that whenever certain quasi-homotheties exist, they can be used to upgrade results from the weak cop game to the strong cop game. An application of these results is given in Corollary \ref{cor:no-homothety-in-hyp-plane}, where we show that the hyperbolic plane does not admit surjective homotheties. A map $h \colon X \to X$ on a metric space $X$ is called a \emph{homothety} if there is some $\mu > 1$ such that $d(h(x),h(y)) = \mu d(x,y)$ for all $x,y \in X$. In analogy to quasi-isometries, we generalize this notion as follows. Note that every homothety $h$ and its powers $h^j$ are already quasi-isometries, but the sequence $h^j$ does not have uniform quasi-isometry constants.

\begin{definition}\label{def:quasi-homothety-seq}
    Let $A\geq 1, B\geq 0$, $\rho = (\rho_j)_j$ a sequence of positive integers with $\rho_j \to \infty$ and $\Gamma,\Delta_j$ graphs for $j\in \mathbb{N}$. A collection of $\iota_j \colon \Delta_j \to \Gamma$ is called a \emph{sequence of $(A\rho + B)$-surjective $(A,B)$-quasi-$\rho$-homotheties} if for all $j\in \mathbb{N}$
\begin{itemize}
    \item [(1)]  $\forall x \in \Gamma , \exists \overline{x} \in \Delta_j \colon
    d_{\Gamma}(\iota_j(\overline{x}),x) \leq A\rho_j + B$.

    \item [(2)] $\forall \overline{x}, \overline{y} \in \Delta_j \colon
    A^{-1}d_{\Delta_j}(\overline{x},\overline{y})- B \leq \frac{1}{\rho_j}d_{\Gamma}(\iota_j(\overline{x}), \iota_j(\overline{y})) \leq A \, d_{\Delta_j}(\overline{x},\overline{y}) + B $.
\end{itemize}
\end{definition}
Note that every surjective homothety $h \colon \Gamma \to \Gamma$ with scaling $\mu >1$ gives rise to a sequence $\iota_j:=h^j \colon \Gamma \to \Gamma$ of $\rho$-surjective $(1,0)$-quasi-$\rho$-homotheties, where $\rho_j = \mu^j$. As for quasi-isometries, quasi-homotheties admit quasi-inverses described in the following Lemma.

\begin{lemma} \label{lem:quasi-homothety-seq}
    Let $\iota_j \colon \Delta_j \to \Gamma$ be a sequence of $(A\rho+B)$-surjective $(A,B)$-quasi-$\rho$-homotheties. Then there are maps $\pi_j \colon \Gamma \to \Delta_j $ such that for all $x,y \in \Gamma$
     $$
    \frac{1}{A} d_{\Delta_j}(\pi_j(x),\pi_j(y)) - (2A+3B) \leq \frac{1}{\rho_j} d_{\Gamma}(x,y) \leq A\, d_{\Delta_j}(\pi_j(x),\pi_j(y)) + (2A + 3B)
    $$
    and for all $x\in \Gamma , \overline{x} \in \Delta_j$
    $$
    d_{\Gamma}(x,\iota_j(\pi_j(x))) \leq A\rho_j+ B  \quad \text{ and } \quad d_{\Delta_j}(\overline{x}, \pi_j(\iota_j(\overline{x}))) \leq AB.
    $$
\end{lemma}
\begin{proof}
    We fix $j$ and for simplicity just write $\Delta=\Delta_j$, $\iota = \iota_j, \rho = \rho_j$ and $\pi = \pi_j$. For $x \in \Gamma$, there exists $\overline{x} \in \Delta$ (without loss of generality let $\iota(\overline{x}) = x$ when $x \in \iota(\Delta)$) such that 
    $$
    d_{\Gamma}(x, \iota(\overline{x})) \leq A\rho +B.
    $$
    Let $\pi(x) := \overline{x}$. Then
    \begin{align*}
    d_{\Gamma}(x,y) & \leq d_{\Gamma}(x, \iota(\overline{x})) + d_{\Gamma}(\iota(\overline{x}), \iota(\overline{y})) + d_{\Gamma}(\iota(\overline{y}), y) \\
    & \leq 2(A\rho + B) + \rho( A\, d_{\Delta}(\overline{x}, \overline{y}) + B ) \\
    & \leq \rho ( A \, d_{\Delta}(\overline{x}, \overline{y}) + 2A + 3B )
    \end{align*}
    and
    \begin{align*}
        d_{\Gamma}(x,y) & \geq  d_{\Gamma}(\iota(\overline{x}), \iota(\overline{y})) - d_{\Gamma}(x, \iota(\overline{x})) - d_{\Gamma}(\iota(\overline{y}), y) \\
        & \geq \rho \left(\frac{1}{A}d_{\Delta}(\overline{x}, \overline{y}) - B\right) - 2(A \rho + B) \\
        & \geq \rho\left(  \frac{1}{A}d_{\Delta}(\overline{x}, \overline{y}) - (2A + 3B) \right) .
    \end{align*}
    It may happen that $\iota$ is not injective, but we still have for $\overline{x} \in \Delta$
    \begin{align*}
        d_{\Delta}(\overline{x} , \pi(\iota(\overline{x} )))
        & \leq  A\left( \frac{1}{\rho} d_{\Gamma}(\iota(\overline{x} ), \iota(\pi(\iota(\overline{x} )) )) + B \right)  \\
        &= A\left( \frac{1}{\rho} d_{\Gamma}(\iota(\overline{x} ), \iota(\overline{x}  )) + B \right) = AB. \qedhere
    \end{align*}
\end{proof}

\begin{remark}
	To elaborate on the analogy to quasi-isometries, consider the following category:
	Objects are sequences of metric spaces and morphisms are uniform sequences of quasi-isometric embeddings, where by uniform we mean that all maps use the same error constants.
	A sequence of quasi-surjective quasi-\(\rho\)-homotheties from \((X, \mu)\) to \((Y, \nu)\) represents an isomorphism from the sequence \((X, \rho_j \mu)_{j\in \mathbb{N}}\) to the constant sequence \((Y, \nu)_{j\in\mathbb{N}}\).
\end{remark}

The following is the most general version of the meta-gaming theorem we provide in this work.

\begin{theorem}[Meta-gaming theorem]\label{thm:meta-gaming}
    Let $\Delta_j$ be graphs with $\wCop(\Delta_j) \geq m \in \mathbb{N}\cup \{\infty\}$ such that for every $\overline{\sigma}, \overline{\rho}$ there exist $\overline{\psi}, \overline{R}$ such that there exist winning strategies for the robber on $\Delta_j$ against $m$ cops, where if the cop picks speed \(\overline{\sigma}\) and reach \(\overline{\rho}\), the robber picks speed $\overline{\psi}$ and reach $\overline{R}$. 
    
    Let $A\geq 1, B\geq 0$ and $\rho_j \to \infty$ a sequence of positive integers. Let $\iota_j \colon \Delta_j \to \Gamma$ be a sequence of $(A\rho+B)$-surjective $(A,B)$-quasi-$\rho$-homotheties. Then $\sCop(\Gamma)  \geq m$.
\end{theorem}

To digest this statement, consider the following, which is \Cref{thm:meta-gaming} if all the \(\Delta_i\) are isometric to some common graph \(\Gamma\).

\begin{corollary}\label{cor:self-meta-gaming}
    Let $\Gamma$ be a graph, $A \geq 1, B\geq 0$ and $\rho_j \to \infty$ a sequence of positive integers. If there is a sequence $\iota_j \colon \Gamma \to \Gamma$ of $(A \rho+B)$-surjective $(A,B)$-quasi-$\rho$-homotheties, then $\sCop(\Gamma) = \wCop(\Gamma)$.
\end{corollary}

For example, for \(\Gamma = \Cay(\Z^2)\), such \(\rho_j\) and fitting homotheties exist. We conclude that \(\sCop(\Z^2) = \wCop(\Z^2)\).
Thus, \Cref{thm:meta-gaming} is really a generalization of \Cref{thm:scop-Zn}. Their proofs use the same basic idea.

\begin{proof}[Proof of \Cref{thm:meta-gaming}]

	Let \(n < m \leq \wCop(\Delta_j)\). For the constants $A,B$ in the statement, let $\overline{\sigma} := 4A^2+3AB$ and $\overline{\rho} = 4A(2A+3B)$. According to the assumptions there exist $\overline{\psi}$ and $\overline{R}$ for which the robber can win on $\Delta_j$ (for any $j$) against $n$ weak cops assuming that the cops started by picking $\overline{\sigma}$ and $\overline{\rho}$.

    We now describe a winning strategy for the robber player against $n$ strong cops on $\Gamma$. The game starts with the cop choosing some speed $\sigma \geq 1$. The robber chooses speed
    $$
    \psi = \sigma  \left( A \overline{\psi} + B\right) .
    $$
    The cop then chooses some reach $\rho$. Without loss of generality, we may assume that $\rho > \sigma$ and that $\rho = \rho_j$ for some $j$ large enough. The robber chooses the radius
    $$
    R = \rho (A\overline{R} + A + 2B).
    $$
    The cop then chooses starting positions $c_i \in \Gamma$ for $i = 1 , \ldots , n$. For the choice of $j$ with $\rho = \rho_j$ above, we fix $\Delta := \Delta_j$ and $\iota := \iota_j \colon \Delta \to \Gamma$. Let $\pi  \colon \Gamma \to \Delta$ be the quasi-projection from Lemma \ref{lem:quasi-homothety-seq}. Definition \ref{def:quasi-homothety-seq} and Lemma \ref{lem:quasi-homothety-seq} are summed up by the following facts for $x,y \in \Gamma, \overline{x}, \overline{y} \in \Delta$.
    \setcounter{equation}{0}
    \begin{align}
        &\forall x \in \Gamma \, \exists \overline{x} \in \Delta \colon d_{\Gamma}(\iota(\overline{x}),x) \leq A\rho + B \label{eq:def-surj} \\
    & \frac{1}{A} d_{\Delta}(\overline{x},\overline{y}) - B \leq \frac{1}{\rho}d_{\Gamma}(\iota(\overline{x}),\iota(\overline{y}))  \leq A \, d_{\Delta}(\overline{x},\overline{y}) + B \label{eq:def-quasi} \\
    &\frac{1}{A}  d_{\Delta}(\pi(x),\pi(y))  -  (2A + 3B) \leq \frac{1}{\rho} d_{\Gamma}(x,y)\leq A\,  d_{\Delta}(\pi(x),\pi(y))  +  (2A + 3B) \label{eq:lem-quasi} \\
    &    d_{\Gamma}(x,\iota(\pi(x))) \leq A\rho+ B \label{eq:lem-surj} \\
    & d_{\Delta}(\overline{x}, \pi(\iota(\overline{x})) ) \leq AB \label{eq:lem-proj}
    \end{align}
    During the game, the robber uses the weak-cop game played on $\Delta$ with parameters $\overline{\sigma}, \overline{\rho}, \overline{\psi} , \overline{R}$ as an oracle for how to play on $\Gamma$. To determine the robber's starting position $r_0$, the robber consults the oracle with cop starting positions $\pi(c_i) \in \Delta$ to obtain a starting position $\overline{r} \in \Delta$ with which the weak-cop game can be won on $\Delta$, in particular $d_\Delta(\overline{r}, \pi(c_i)) > \overline{\rho} + \overline{\sigma}$. The robber then chooses the actual starting position $r := \iota(\overline{r}) \in \Gamma$. By (\ref{eq:lem-quasi}) and (\ref{eq:lem-proj}),
    \begin{align*}
        d_{\Gamma}(\iota(\overline{r}),c_i) & \geq \rho \left( \frac{1}{A} \left( d_\Delta (\pi(\iota(\overline{r})), \pi(c_i) )- (2A + 3B) \right)\right) \\
        &\geq \rho \left( \frac{1}{A} ( d_\Delta(\overline{r}, \pi(c_i)) - d_\Delta(\overline{r} , \pi(\iota(\overline{r})) ) )-(2A+ 3B) \right) \\
        &> \rho \left( \frac{1}{A} (\overline{\rho} + \overline{\sigma} - AB ) -(2A + 3B) \right) \\
        &=  \rho\left( 4(2A + 3B) + \frac{1}{A}\overline{\sigma} - B - (2A + 3B) \right) \\
        & \geq \rho(2(2A + 3B)) \geq 2\rho \geq \rho + \sigma ,
    \end{align*}
    the robber is not caught at the starting position $r$, even after the cops move once.

    The game now proceeds in meta-stages. At the beginning of a meta-stage, the robber is placed at $r_t = \iota(\overline{r}_t)$ for some $\overline{r}_t \in \Delta$. The cops' positions $c_i$ are projected to $\Delta$, and the oracle suggests a path $ \overline{r}_t = \overline{p}_1 , \overline{p}_2, \ldots , \overline{p}_k \in \Delta$ along which the robber follows in a winning strategy on $\Delta$. The actual robber will now spend the next $\lambda := \lceil \rho/\sigma \rceil $ many turns following geodesic paths from $p_j := \iota(\overline{p}_j)$ to $p_{j+1}$ for $j \in \{ 1, \ldots , k-1\}$ to finally reach $p_k \in \iota(\Delta)$. During the meta-stage, the cops can also do $\lambda$ many moves. The meta-stage strategy then repeats. To prove that this gives a winning strategy for the robber it remains to show that
    \begin{itemize}
        \item [(a)] the robber can actually move as described from $p_1$ to $p_k$ in $\lambda$ many moves,
        \item [(b)] the movement of the cops for $\lambda$ many turns corresponds to a valid single move when projected to $\Delta$,
        \item [(c)] the robber is never caught during a meta-stage, and
        \item [(d)] after any number of turns, the robber eventually returns to the ball of radius $R$ around the base point.
    \end{itemize}

    (a) We have $d_\Delta(\overline{p}_1, \overline{p}_k) \leq \overline{\psi}$. Using (\ref{eq:def-quasi}) and $1 \leq \rho/\sigma \leq \lceil \rho / \sigma \rceil = \lambda$,
\begin{align*}
    d_{\Gamma}(p_1, p_k) & \leq \rho\left(A \, d_\Delta(\overline{p}_1, \overline{p}_k) + B \right)
    \leq \frac{\rho}{\sigma} \sigma \left(A \overline{\psi} + B\right)
    \leq \lambda \sigma (A \overline{\psi} + B) \leq \lambda \psi
\end{align*}
shows that the robber can complete the movement from $p_1$ to $p_k$ within $\lambda$ many turns.

\begin{figure}
  \centering
  \includegraphics[width=1\textwidth]{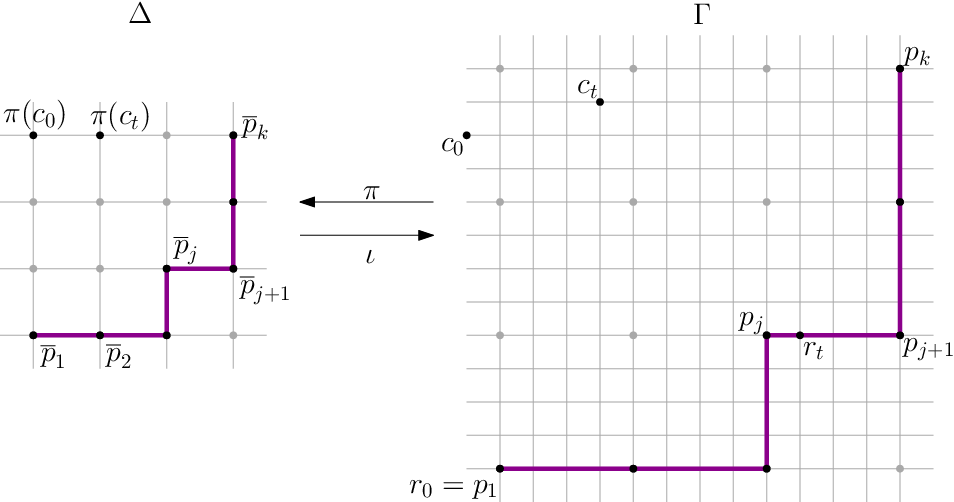} 
  \caption{At the beginning of a meta-stage, the cop's position $c_0$ is provided to the oracle $\Delta$ via $\pi$. In the oracle there is a winning strategy by following along the path $\overline{p}_1, \ldots , \overline{p}_k$, which the robber implements in $\Gamma$, by following along $p_1 := \iota(\overline{p}_1), \ldots , p_k$. }
  \label{fig:meta_game}
\end{figure}
(b) If the starting position of a cop is $c_0\in \Gamma$, then their position $c_t \in \Gamma$ after $\lambda$ many turns satisfies $d_{\Gamma}(c_0,c_t) \leq \lambda \sigma $. Projected to $\Delta$, we use (\ref{eq:lem-quasi}) and $\lambda \sigma = \lceil \rho/\sigma \rceil \sigma \leq (\rho/\sigma + 1) \sigma = \rho + \sigma \leq 2\rho$ to show
\begin{align*}
    d_\Delta(\pi(c_0), \pi(c_t))
    & \leq A \left( \frac{1}{\rho} d_{\Gamma}(c_0,c_t) + 2A + 3B\right)
    \\
    &\leq A \left( 2 + 2A +3B \right) \leq A(4A +3B) = \overline{\sigma}
\end{align*}
which confirms that the projected movement of the cops corresponds to a valid movement in $\Delta$.

(c) Let $c_t$ be the position of a cop and $r_t$ the position of the robber at any time (possibly even different times) during a meta-stage. Note that in the projected game, the cops $\pi(c_t)$ are never too close to any of the vertices $\overline{p}_1, \ldots , \overline{p}_k$ of the path of the robber, namely
$$
d_\Delta(\pi(c_t), \overline{p}_i) > \overline{\rho}
$$
for every $i \in \{1, \ldots, k\}$. Since the robber in $\Gamma$ moves on geodesics along the $p_i$, there is a $j \in \{1, \ldots, k\}$ such that $d_{\Gamma}(r_t, p_j) \leq d_{\Gamma}(p_j, p_{j+1})$. The situation is illustrated in Figure \ref{fig:meta_game}. Using (\ref{eq:def-quasi}), we get an estimate
\begin{align*}
    d_{\Gamma}(r_t, p_j)
    &\leq d_{\Gamma}(\iota(\overline{p}_j), \iota(\overline{p}_{j+1}))
     \leq \rho ( A d_\Delta(\overline{p}_j, \overline{p}_{j+1}) + B ) = \rho(A + B)
\end{align*}
on the distance from the robber to some $p_j$. We then use the triangle inequality, the above estimates, (\ref{eq:def-quasi}), (\ref{eq:lem-surj}) and $\overline{\rho} = 4A(2A + 3B)$ to obtain
\begin{align*}
    d_{\Gamma}(r_t, c_t) & \geq d_{\Gamma}(\iota (\pi (c_t)) , \iota(\overline{p}_j)) - d_{\Gamma}(r_t , p_j) - d_{\Gamma}(\iota(\pi(c_t)), c_t) \\
    & \geq \rho \left(\frac{1}{A} d_\Delta(\pi(c_t), \overline{p}_j) - B \right) - \rho (A + B) - (A\rho + B) \\
    & > \rho \left( \frac{1}{A} \overline{\rho} - B \right) - 2\rho(A+B) \\
    & = \rho ( 4(2A+3B) -B -2A -2B) \geq 2\rho \geq \rho + \sigma
\end{align*}
which means that $r_t$ is always at least $\rho$ far away from $c_t$, even after one more movement by a distance of at most $\sigma$ by the cop. This means that the robber is never caught.

(d) Let $v$ be the base point in $\Gamma$ and $\pi(v)$ the base point in $\Delta$.
For any number of turns $k$, there is a number $k' \geq k$ such that the robber is in the ball of radius $\overline{R}$ around $\pi(v)$, that is $d_\Delta(\pi(v), \overline{r}) \leq \overline{R}$, where $\overline{r}$ is the position of the robber at the beginning of turn $k'$. In $\Gamma$, the robber will eventually visit $\iota(\overline{r})$, which lies in the ball by the following estimate. We use the triangle inequality, (\ref{eq:lem-surj}) and (\ref{eq:def-quasi}) to obtain
\begin{align*}
    d_{\Gamma}(v,\iota(\overline{r})) & \leq d_{\Gamma}(v, \iota(\pi(v))) + d_{\Gamma}( \iota(\pi(v)), \iota(\overline{r}))  \\
    &\leq ( A\rho + B) + \rho(A \, d_\Delta(\pi(v),\overline{r}) + B) \\
    &
    \leq  \rho (A + B + A\overline{R} + B)= R.
\end{align*}

The above calculations show that the robber can implement the strategy from the weak cop game on $\Delta$ using meta-stages to win the strong cop game on $\Gamma$. Hence $\operatorname{sCop}(\Gamma) > n$ for every $n< m = \wCop(\Delta)$, so $\sCop(\Gamma) \geq m$.
\end{proof}

\section{Applications}\label{sec:applications}

In this section, we give some applications of the Meta-gaming Theorem. First, we prove that certain lamplighter groups have infinite strong cop numbers, using the full power of Theorem \ref{thm:meta-gaming}. We then deduce that hyperbolic groups that are not virtually free and the hyperbolic plane do not admit self-homotheties. In Subsection \ref{sec:bs} we use the idea of meta-stages again to calculate the strong cop number of Baumslag-Solitar groups (but this time without homotheties). Finally, we apply Theorem \ref{thm:scop-Zn} to Thompson's group $F$.

\subsection{Lamplighter groups}

Cornect and Martínez-Pedroza proved in \cite[Theorem 1.5]{CoMa24arxiv} that restricted wreath products $L \wr H$ of finitely generated groups where $L$ non-trivial and $H$ infinite have infinite cop number. We reprove their result for $L$ finite in Lemma \ref{lem:wcop(lamplighter)}, keeping track of the constants to be able to use \Cref{thm:meta-gaming}, the meta-gaming theorem, to upgrade the result to the strong cop number when $L$ finite and $H=\mathbb{Z}$, see Theorem \ref{thm:scop-lamplighter}. We start by recalling the lamplighter group.

For finitely generated groups \(L\) and \(H\), we denote their \emph{restricted wreath product} by \(L \wr H \coloneq L^{\oplus H} \rtimes H\), where \(H\) acts by index-shifts.
We denote elements of \(L^{\oplus H}\) as formal sums \(\sum_{h \in H} l_h~h\). This provides a natural embedding \(L \into L^{\oplus H}\) by sending \(l \in L\) to \(l\cdot 1\).
Note that if \(T\) is a generating set for \(H\), then \(S \coloneq L \times (T \cup \{1\})\) is a generating set for \(L \wr H\).
We think of this group as a street of shape \(H\), where at every point there is a lamp with possible states \(L\) and at some point \(p \in H\) there is a lamplighter.
The generating set \(S\) corresponds to the lamplighter switching the lamp at their current position to any state and then moving along the street using some generator \(t \in T\).
Hence we also call \(L \wr H\) a \emph{lamplighter group}.

\begin{lemma}\label{lem:wcop(lamplighter)}
	Let \(L\) be a non-trivial finite group and \(H\) a finitely generated infinite group and set \(G \coloneq L \wr H\).
	Then \(\wCop(G) = \infty\). Moreover, given the cop's speed $\sigma$ and reach $\rho$, the robber can choose speed $\psi$ and radius $R$ depending only on $H$ and not on $L$ to win the weak-cop game.
\end{lemma}
\begin{proof}
	Let \(T\) be a finite generating set for \(H\) and consider the finite generating set \(S \coloneq L \times (T \cup \{1\})\) for $G$.
	Let \(n \in \N\). We provide a winning strategy for the robber against \(n\) weak cops on \( \Gamma := \Cay(G,S)\).

	Note that an element of \(\Gamma\) corresponds to a street with lamps and a lamplighter. Having the lamplighter move around switching lamps corresponds to a path in \(\Gamma\). Forgetting the states of the lamps gives paths on $\Delta := \Cay(H,T)$.

	Suppose that the cops picked speed \(\sigma\) and reach \(\rho\). 
	We choose \(2n\) pairwise distinct elements \(a_i, b_i \in \Delta\) such that \(d_\Delta(a_i, b_i) > \sigma + \rho\). This ensures that for any \(g \in \Gamma\) and any element of \(h \in B_{\sigma + \rho}(g)\), the states of the lamps at least at one of \(a_i\) and \(b_i\) coincide between \(g\) and \(h\).
	That is, a cop positioned at \(g\) can not, within one move, reach a position where they can see a robber at an element \(h\) which has different lamp-states in both \(a_i\) and \(b_i\).

	Defining $a_0:= a_n$ and $b_0:= b_n$ for convenience, the robber then chooses speed and radius
    $$
    \psi := R:= 2 d_\Delta(a_n, b_n) +  \sum_{i=1}^{n} d_\Delta(a_{i-1},a_i) + d_\Delta(b_{i-1}, b_i)
    $$
    which corresponds to the length of a loop $\gamma$ in $\Delta$ visiting all the points $a_i, b_i$. Note that the numbers $\psi$ and $R$ only depend on $\Delta$ but not on $L$.

    Next, the cops choose their starting positions.
    Let $c_i(p)\in L$ be the lamp-state of the $i$'th cop's lamp at position $p \in \Delta$. The robber chooses \(1 \neq l_0 \in L\) and sets \(\widebar 1 \coloneq l_0\) and \(\widebar l \coloneq 1\) for every \(l \in L \setminus 1\). The robber then chooses a starting position such that the robber's lamps have state $r(a_i)=\widebar{c_i(a_i)}$ and $r(b_i)=\widebar{c_i(b_i)}$, and the lamplighter-position is $a_1$.
	In particular, this is a position where the robber's lamp-states at \(a_i\) and \(b_i\) differ from those of the \(i\)'th cop.
	Hence by construction, the robber cannot be caught in the cop player's first move.

	But, after a cop move, it might be that for some \(i\), \(c_i(a_i)\) equals \(r(a_i)\), the state of the robber's lamp at \(a_i\) and the same goes for \(b_i\). But for any given \(i\), only one of \(c_i(a_i)\) and \(c_i(b_i)\) can match the robber.

	The robber's move is to have the lamplighter move along the loop $\gamma$ traversing all the \(a_i\) and \(b_i\). For every \(a_i\) where \(c_i(a_i) = r(a_i)\), the lamplighter switches the lamp to \(\widebar{c(a_i)}\) and analogously for \(b_i\). As a generator in \(S\) corresponds to switching a lamp and then moving the lamplighter along some \(t \in T\), the robber's move is of length at most \(\psi\).
	By construction, no point on this path can be seen by any of the cops, so this gives a valid path for the robber player. To see this, note that from an individual cop's perspective, any lamp configuration along the path might as well have been the robber's starting position, and the cops are far enough away from the robber that changing the robber's lamplighter-position does not get them close enough.
	Moreover, at no point does the robber leave the ball of radius $R$ around their starting point. Repeating this strategy, the robber is never caught and never leaves the ball, thus winning the game.
\end{proof}

\begin{remark}
    In the above proof, the points $a_i, b_i$ could be chosen in such a way that for some base point $v\in \Delta$ we have $d(v,a_i)=i$ and $d(v,b_i)=\sigma+\rho+i+1$. Then
    $$
    \psi \leq \sum_{i=1}^n 2d(v,a_i) + 2d(v,b_i) \leq 2\sum_{i=1}^n i + (\sigma + \rho +i + 1) = 2n(n+1 + \sigma +\rho + 1),
    $$
    so $\psi$ does not even depend on $\Delta$, but only on $n, \sigma$ and $\rho$.
\end{remark}

In the following, we restrict to $H=\mathbb{Z}$ and $L$ finite. Elements of $G = L\wr \mathbb{Z}$ are then of the form $(a,n)$, where $a = (\ldots , a_0, a_1, \ldots)$ with $a_k \in L$ for $k,n \in \mathbb{Z}$.

\begin{theorem}\label{thm:scop-lamplighter}
	Let \(L\) be a non-trivial finite group and set \(G \coloneq L \wr \Z\).
	Then \(\sCop(G) = \infty\).
\end{theorem}

\begin{proof}
	We use the generating set \(L \times \{t\}\), where \(t\) is a generator of \(\Z\) and let $\Gamma := \Cay(G,L\times \{t\})$. The goal is to apply \Cref{thm:meta-gaming}, so we have to provide a sequence \(\Delta_j\) and a sequence of quasi-homotheties from \(\Delta_j\) to \(\Gamma\), see Definition \ref{def:quasi-homothety-seq}.

	We set \(H_j \coloneq L^j \wr \Z\) and identify it with its Cayley graph $\Delta_j := \operatorname{Cay}(H_j, L^j \times \{s\})$ with respect to the generating set \(L^j \times \{s\}\), where $s$ is a generator of $\mathbb{Z}$.
	We further set \(\rho_j \coloneq j\) and proceed to show that there exists a sequence of quasi-surjective quasi-\(\rho\)-homotheties \(\Delta_j \to \Gamma\).
	The maps are given by $\iota_j \colon \Delta_j \to \Gamma,$ \[
		\iota_j(a, n) \coloneq ( b, j \cdot n)
	\] where \(b_{j\cdot k + i} = a_k^{(i)}\) for $a_k = (a_k^{(0)}, \ldots , a_k^{(j-1)}) \in L^j$ and any \(k \in \Z\) and \(i \in \{0 , \dots , j-1\}\). See an illustration of this in \Cref{fig-lamplighter-homothety}. The constants that appear in Definition \ref{def:quasi-homothety-seq} are $A=1, B=2$.

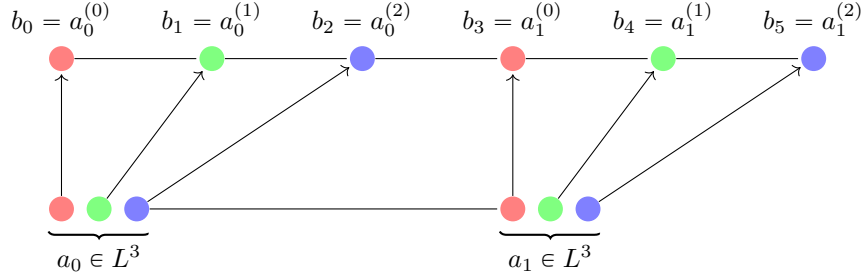
\begin{figure}
	\centering
	\begin{tikzpicture}[shorten >=1pt,on grid,auto]

	\node[circle,label={\(b_0=a^{(0)}_0\)},fill,red!50]   (G00) {};
	\node[circle,label={\(b_1=a^{(1)}_0\)},fill,green!50] (G01) [right=2 of G00] {};
	\node[circle,label={\(b_2=a^{(2)}_0\)},fill,blue!50]  (G02) [right=2 of G01] {};
	\node[circle,label={\(b_3=a^{(0)}_1\)},fill,red!50]   (G10) [right=2 of G02] {};
	\node[circle,label={\(b_4=a^{(1)}_1\)},fill,green!50] (G11) [right=2 of G10] {};
	\node[circle,label={\(b_5=a^{(2)}_1\)},fill,blue!50]  (G12) [right=2 of G11] {};

	\node[circle,fill,red!50]   (H00) [below=2 of G00] {};
	\node[circle,fill,green!50] (H01) [right=.5 of H00] {};
	\node[circle,fill,blue!50]  (H02) [right=.5 of H01] {};

	\node[circle,fill,red!50]   (H10) [below=2 of G10] {};
	\node[circle,fill,green!50] (H11) [right=.5 of H10] {};
	\node[circle,fill,blue!50]  (H12) [right=.5 of H11] {};

	\path[->]
		(H00) edge (G00)
		(H01) edge (G01)
		(H02) edge (G02)
		(H10) edge (G10)
		(H11) edge (G11)
		(H12) edge (G12)
	;
	\path[-]
		(G00) edge (G01)
		(G01) edge (G02)
		(G02) edge (G10)
		(G10) edge (G11)
		(G11) edge (G12)

		(H02) edge (H10)
	;
	\draw [
		thick,
		decoration={
			brace,
			mirror,
			raise=0.3cm
		},
	decorate
	] (H00.west) -- (H02.east)
	node [pos=0.5,anchor=north,yshift=-0.35cm] {\(a_0 \in L^3\)};
	\draw [
		thick,
		decoration={
			brace,
			mirror,
			raise=0.3cm
		},
	decorate
	] (H10.west) -- (H12.east)
	node [pos=0.5,anchor=north,yshift=-0.35cm] {\(a_1 \in L^3\)};

\end{tikzpicture}
	\caption{The map \(\iota_3 \colon L^3 \wr \Z \to L \wr \Z\) of a quasi-homothety.}
	\label{fig-lamplighter-homothety}
\end{figure}

	To prove quasi-surjectivity, note that \((b, n) \in \iota_j(\Delta_j)\) if and only if \(n\) is a multiple of \(j = \rho_j\). If \(n'\) is the largest multiple of \(j\) such that \(n' \leq n\), we have \(d_\Gamma \left( (b, n), (b, n') \right) = n - n' \leq \rho_j\) and \((b, n') \in \iota_j(\Delta_j)\).

	For the quasi-homothety-part let \((a, n), (a', n') \in \Delta_j\).
    If $n\neq n'$ or $a_k \neq a_k'$ for some $k\neq n$, the distance between these two elements in \(\Delta_j\) is the length of the shortest lamplighter-path in \(\Z\) from \(n\) to \(n'\) that traverses all indices $k$ where \(a_k \neq a_k'\).
	Similarly, the distance \(d_\Gamma\left(\iota_j(a, n), \iota_j(a', n')\right)\) is the length of the shortest lamplighter-path from \(jn\) to \(jn'\) that traverses all indices $k$ where \(b_k \neq b_k'\). If $n=n'$ and for all $k\neq n$ we have $a_{k} = a'_k$, then the distance from $(a,n)$ to $(a',n')$ is $2$ and the lamplighter has to travel along $tt^{-1}$, while the distance from $\iota_j(a,n)$ to $\iota_j(a',n')$ is between $2$ and $2j$.

	Up to replacing $t$ with $t^{-1}$, we may assume without loss of generality that $n \leq n'$.
    In both cases above, the lamplighter-path is of the form $t^{-r}t^{\ell}t^{-m}$ for some \(r, \ell, m \in \mathbb{Z}_{\geq 0}\) such that \(-r+\ell-m = n'-n \).
    If for some \(k \in \Z\) we have \(b_k \neq b'_k\), then \(a_{\floor{k/j}} \neq a'_{\floor{k/j}}\). Hence there is a (not necessarily shortest) lamplighter-path of the form $s^{-jr}s^{j(\ell+1)} s^{-j(m+1)}$ that realizes a path from $\iota_j(a,n)$ to $\iota_j(a',n')$ in $\Gamma$. It has length at most \(j(r+\ell+m+2)\), so
    $$
    d_\Gamma(\iota_j(a,n), \iota_j(a',n')) \leq j(d_{\Delta_j}((a,n),(a',n')) + 2 ).
    $$
    On the other hand, every lamplighter-path corresponding to a path from \(\iota_j(a, n)\) to \(\iota_j(a', n')\) has to traverse at least the points $jn$, \(j(n - r)\), \(j(n - r + \ell)\) and \(jn'\) in order, so has length at least \(j(r+\ell+m)\).
	Hence \begin{align*}
		d_{\Delta_j}\left((a, n), (a', n')\right) &= r+\ell+m \leq \frac 1j d_\Gamma\left(\iota_j(a, n), \iota_j(a', n')\right) \leq d_{\Delta_j}\left((a, n), (a', n')\right)
	\end{align*} and \(\iota_j\) is indeed a sequence of $(1,2)$-quasi-surjective $\rho$-homotheties. Since the speed $\psi$ and radius $R$ in the weak cop game only depends on the speed $\sigma$, reach $\rho$ and the group $\mathbb{Z}$, but not on $\Delta_j$ by Lemma \ref{lem:wcop(lamplighter)}, we can apply Theorem \ref{thm:meta-gaming} to conclude $\sCop(G)= \infty$.
\end{proof}

\newpage 

\subsection{Non-existence of quasi-homotheties}
\begin{corollary}\label{cor:no-homothety-in-hyp-groups}
    If $\Gamma$ is the Cayley graph of a hyperbolic group $G$ that is not virtually free, then there does not exist a sequence $\iota_j \colon \Gamma \to \Gamma$ of $(A \rho+B)$-surjective $(A,B)$-quasi-$\rho$-homotheties, for any $A\geq 1, B\geq 0, \rho_j \to \infty$. For $A=1, B=0$, there does not exist a sequence of $\rho$-surjective $\rho$-homotheties for $\rho_j \to \infty$.
\end{corollary}
\begin{proof}
By Theorems \ref{thm:scop-hyperbolic} and \ref{thm:wcop-virtfree}, for such graphs $\Gamma$, we have $\wCop(\Gamma) > 1$, but $\sCop(\Gamma) = 1$. The existence of a self-quasi-homothety would contradict Corollary \ref{cor:self-meta-gaming}.
\end{proof}
\begin{corollary}\label{cor:no-homothety-in-hyp-plane}
    The hyperbolic plane $\mathbb{H}^2$ does not admit a surjective homothety $h \colon \mathbb{H}^2 \to \mathbb{H}^2$ with $\mu > 1$.
\end{corollary}
\begin{proof}
It is well known that the hyperbolic plane is quasi-isometric to the Cayley graph of a cocompact lattice $G < \operatorname{Isom}(\mathbb{H}^2)$, and that all such lattices are hence one-ended, non-amenable and hyperbolic. Let $\Gamma$ be a Cayley graph of $G$ and $f \colon \Gamma \to \mathbb{H}^2, g \colon \mathbb{H}^2 \to \Gamma$ be $(C,D)$-quasi-isometries for $C\geq 1, D \geq 0$. We will now show that if $h \colon \mathbb{H}^2 \to \mathbb{H}^2$ is a surjective $\mu$-homothety for $\mu>1$, then $\iota_j := g \circ h^j \circ f \colon \Gamma \to \Gamma$ is a sequence of $(A\rho+B)$-surjective $(A,B)$-quasi-$\rho$-homotheties for $\rho_j = \mu^j$, $A= C^2D + 1$ and $B= CD+D$. But by Corollary \ref{cor:no-homothety-in-hyp-groups}, this is impossible, so no surjective $\mu$-homothety of the hyperbolic plane exists.

The $(C,D)$-homotheties $f\colon \Gamma \to \mathbb{H}^2$ and $g \colon \mathbb{H}^2 \to \Gamma$ satisfy for $\overline{x}, \overline{y} \in \Gamma$ and $x,y \in \mathbb{H}^2$
\setcounter{equation}{0}
\begin{align}
& \frac{1}{C}d_{\Gamma}(\overline{x},\overline{y}) - D  \leq d_{\mathbb{H}^2}(f(\overline{x}), f(\overline{y}) )  \leq C \, d_{\Gamma}(\overline{x},\overline{y})  + D \label{eq:hyp-qi-f}\\
		&\forall y \in \mathbb{H}^2 \ \exists \overline{x} \in {\Gamma} \colon d_{\mathbb{H}^2}(f(\overline{x}),y) < D \label{eq:hyp-surj-f}\\
    &\frac{1}{C}d_{\mathbb{H}^2}(x,y) - D  \leq d_{\Gamma}(g(x),g(y) )  \leq C \, d_{\mathbb{H}^2}(x,y) + D \label{eq:hyp-qi-g} \\
		&\forall \overline x \in {\Gamma} \ \exists  y \in \mathbb{H}^2   \colon d_{\Gamma}(g(y),\overline x) < D . \label{eq:hyp-surj-g}
\end{align}
We check the conditions in \Cref{def:quasi-homothety-seq}. The situation is illustrated in the following commutative diagram.
\begin{center}
\begin{tikzcd}
\overline{y} \in\mathbb{H}^2 \arrow[r, "h^j"]         & y \in \mathbb{H}^2 \arrow[d, "g"] \\
\overline{x}\in {\Gamma} \arrow[u, "f"] \arrow[r, "\iota_j"] & x \in {\Gamma}
\end{tikzcd}
\end{center}
Let $x \in {\Gamma}$. By (\ref{eq:hyp-surj-g}) there exists $y \in \mathbb{H}^2$ with $d_{\Gamma}(g(y),x) < D$. Since $h$ and $h^j$ are surjective there is some $\overline{y}\in \mathbb{H}^2$ with $h^j(\overline{y}) = y$. By (\ref{eq:hyp-qi-f}) there is $\overline{x} \in {\Gamma}$ with $d_{\mathbb{H}^2}(f(\overline{x}), \overline{y}) < D$. We have found $\overline{x }\in {\Gamma}$ with
\begin{align*}
    d_{\Gamma}(\iota_j(\overline{x}) ,x ) &\leq d_{\Gamma}(\iota_j(\overline{x}), g(y)) + d_{\Gamma}(g(y),x) < C d_{\mathbb{H}^2}(h^j(f(\overline{x})), y) + D + D \\
    &\leq C \mu^j d_{\mathbb{H}^2}(f(\overline{x}), \overline{y}) + 2D < C \mu^j D + 2D \\
    &\leq (C^2 D+1)\mu^j + (CD+D) =  A \rho_j + B
\end{align*}
where we additionally used the triangle inequality and (\ref{eq:hyp-surj-g}). Now if $\overline{x}, \overline{y} \in {\Gamma}$, we use (\ref{eq:hyp-qi-g}) and (\ref{eq:hyp-qi-f}) to obtain
\begin{align*}
    d_{\Gamma}(\iota_j(\overline{x}), \iota_j(\overline{y}))
    &\leq  C d_{\mathbb{H}^2}(h^j(f(\overline{x})), h^j(f(\overline{y}))) + D
    =    C \mu^j d_{\mathbb{H}^2}(f(\overline{x}), f(\overline{y})) + D \\
    & \leq C \mu^j (C d_{\Gamma}(\overline{x}, \overline{y}) + D) + D
     \leq (C^2D+1) \mu^j d_{\Gamma}(\overline{x}, \overline{y}) + \mu^j (CD + D)
\end{align*}
and
\begin{align*}
    d_{\Gamma}(\iota_j(\overline{x}), \iota_j(\overline{y}))
    &\geq \frac{1}{C} d_{\mathbb{H}^2}(h^j(f(\overline{x})), h^j(f(\overline{y}))) - D
     \geq \frac{1}{C} \mu^j d_{\mathbb{H}^2}(f(\overline{x}), f(\overline{y})) - D \\
    & \geq \frac{\mu^j}{C}  \left( \frac{1}{C} d_{\Gamma}(\overline{x},\overline{y}) - D \right) - D
		\geq   \frac{1}{C^2D + 1}\mu^j d_{\Gamma}(\overline{x},\overline{y}) -\mu^j (CD + D)
\end{align*}
whence
$$
A^{-1}d_{\Gamma}(\overline{x},\overline{y})- B \leq \frac{1}{\rho_j}d_{\Gamma}(\iota_j(\overline{x}), \iota_j(\overline{y})) \leq A \, d_{\Gamma}(\overline{x},\overline{y}) + B,
$$
concluding the proof that $\iota_j$ is a sequence of $(A\rho+B)$-surjective $(A,B)$-quasi-$\rho$-homotheties.
\end{proof}

\subsection{Baumslag-Solitar groups $\operatorname{BS}(1,n)$}\label{sec:bs}

\begin{figure}[h]
  \centering
  \includegraphics[width=0.3\textwidth]{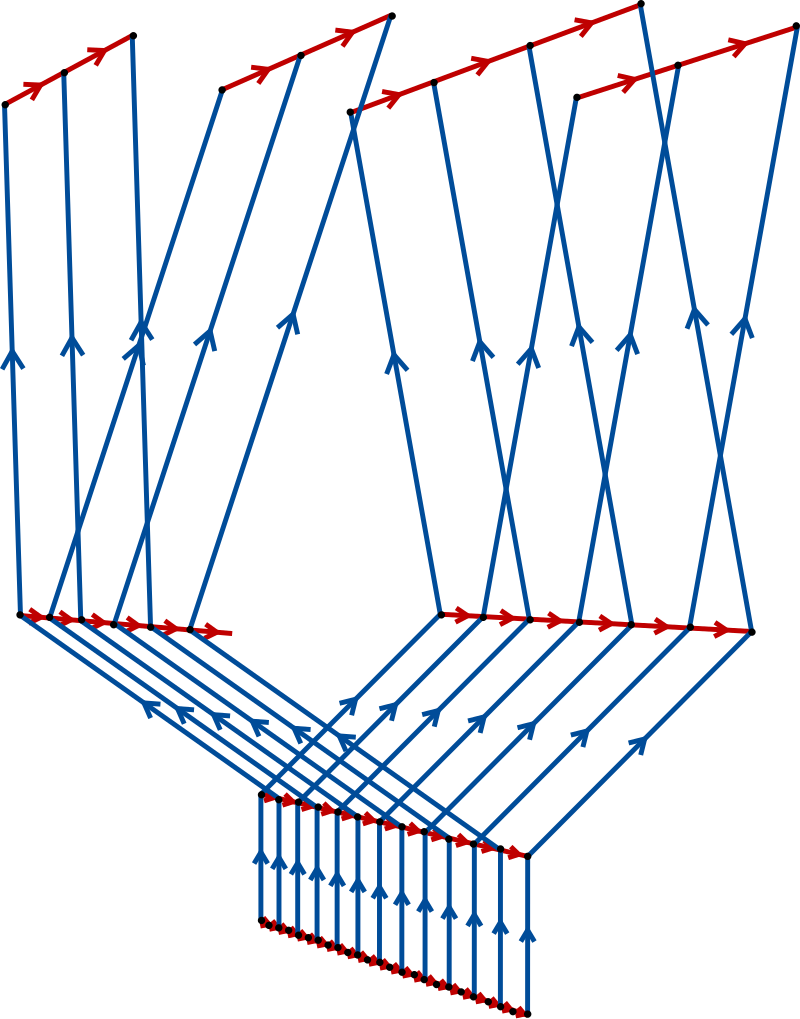}
  \caption{Parts of the Cayley graph for the Baumslag-Solitar group $\operatorname{BS}(1,2) = \langle a,t \colon tat^{-1} = a^2 \rangle$ with generating set $\{a,t\}$.}
  \label{fig:bs}
\end{figure}

The \emph{solvable Baumslag-Solitar groups}
$$
\operatorname{BS}(1,n) := \left\langle a,t \colon tat^{-1} = a^n\right\rangle
$$
for $n \geq 1$ are a well studied class of groups that include $\Z^2$ for $n=1$. A part of the Cayley graph of $\operatorname{BS}(1,2)$ with respect to the standard generating set $\{a,t\}$ is depicted in Figure \ref{fig:bs}. When $n \geq 2$, the Baumslag-Solitar groups are known to have some properties of non-positive curvature, such as exponential growth \cite{CEG94} and exponential divergence \cite{sisto2012arxiv}, but they are not Gromov-hyperbolic. We show that $\sCop(\operatorname{BS}(1,n)) = \infty$.

\begin{theorem}\label{thm:scop-bs}
For all $m \geq 1$, $\sCop(\operatorname{BS}(1,m)) = \infty$.
\end{theorem}
\begin{proof}
We note that $\operatorname{BS}(1,1) \cong \Z^2$, so the results holds for $m=1$ by Theorem \ref{thm:scop-Zn}. We will show below that $\sCop(\operatorname{BS}(1,n+1)) > n$ for all $n\geq 1$. This implies that whenever $n < m$, then $\sCop(\operatorname{BS}(1,m)) > n$, as the robber player can just pretend that the missing $m-1-n$ cops are placed at $1$ and stay there during the whole game. Now for arbitrary $m \geq 2$ and $n\geq 1$ we can find $k\geq 1$ such that $n < m^k$. Farb and Mosher showed in \cite{FaMo98} that $\operatorname{BS}(1,m)$ is quasi-isometric to $\operatorname{BS}(1,m^k)$. Since the cop numbers are invariant under quasi-isometry \cite[Corollary G]{LMR23}, we then have $\sCop(\operatorname{BS}(1,m)) = \sCop(\operatorname{BS}(1,m^k)) > n$. This is true for all $n$, so $\sCop(\operatorname{BS}(1,m)) = \infty$ for all $m \geq 2$ as well.
It thus suffices to show that there is a winning strategy for the robber player against $n$ cops on the Cayley graph of $\operatorname{BS}(1,n+1)$, which we will do in the following.

    We assume the base point is at the neutral element. First, the cops choose a speed $\sigma$, without loss of generality we may assume $\sigma \geq n+1$. The robber chooses $\psi = 34\sigma$. The cops then choose a reach $\rho$, without loss of generality we may assume that $3\sigma < \rho$. Finally, the robber chooses $R = n(n+1)^{8\rho} + 8\rho + n$. Let
    $$
    S_0 := \left\{ t^m a^{N}  \in \operatorname{BS}(1,n+1) \colon 0 \leq m \leq 8\rho, 0 \leq N \leq n(n+1)^{8\rho - m}  \right\}
    $$
    be a rectangular region of the first sheet that has height $8\rho$, length $n$ at the top and length $n(n+1)^{8\rho}$ at the bottom. We call $A = \{a^N \in \operatorname{BS}(1,n+1) \colon N \in \Z\}$ the \emph{axis}. The translates $S_i := a^i S_0$ for $1 \leq i \leq n $ are copies of the rectangular region, but they are part of other sheets, in fact $S_i \cap S_j$ is contained in the axis $A$ for all $i\neq j$. Any path from one sheet to another has to pass through the axis. At the beginning of the game, the cops are placed somewhere and since there are $n+1$ upper parts of distinct sheets, there is at least one (with index $i_0$) that does not contain a cop and the robber chooses to start at $a^{i_0}t^{8\rho}$. The robber is not captured at the beginning and during the first movement of the cops. The robber's strategy now is the following. The robber stays put until one of the cops reaches at least height $\rho/2$ in the upper part of the sheet in which the robber is placed. As soon as that happens, the robber picks one of the other $n$ sheets (say with index $j$) that does not have any cops in the upper part and walks to the point $a^jt^{8\rho}$ along one of $n+1$ possible paths described below. Walking along this path may take more than one turn, but the robber just keeps going until they reach their goal, we call the whole movement one \emph{meta-step}. Then the procedure repeats.

\begin{figure}
  \centering
  \includegraphics[width=0.6\textwidth]{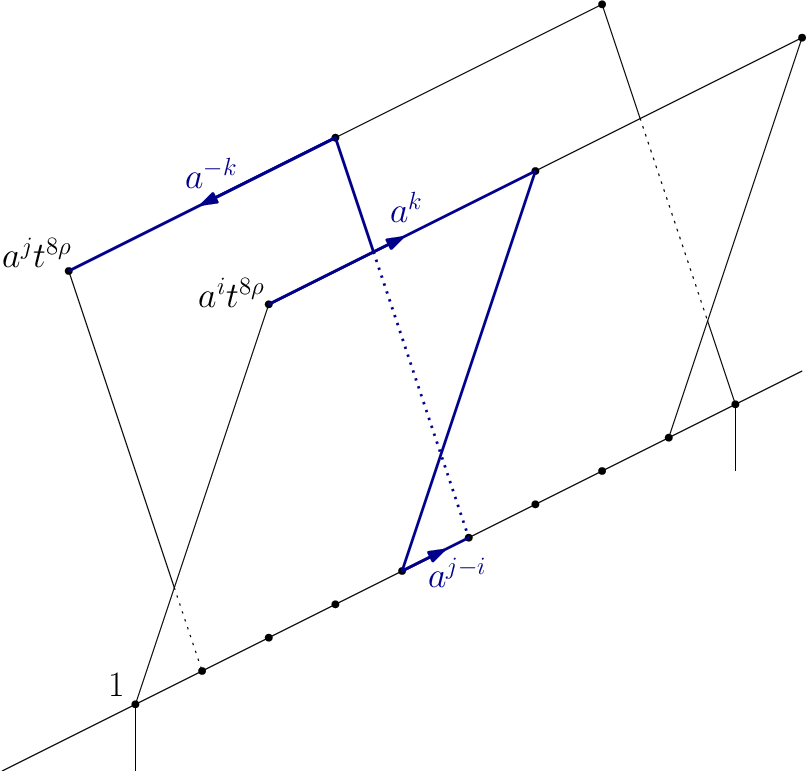} 
  \caption{The path with index $k$ from sheet $S_i$ to sheet $S_j$ passes through the axis $A$. The situation is illustrated for $n=2$, $k=1$, $i=0$, $j=1$, $\rho=1/4$.}
  \label{fig:bs_paths}
\end{figure}

    The $n+1$ possible paths going from $a^it^{8\rho}$ to $a^jt^{8\rho}$, indexed by $k \in \{0,1,\ldots n\}$, are defined by following the sequence of steps $a^kt^{-8\rho}a^{j-i}t^{8\rho}a^{-k}$, see Figure \ref{fig:bs_paths}. The length of one of these paths is at most $3k+16\rho \leq 3n + 16\rho \leq 17\rho$. During the time it takes the robber to move distance $17\rho$, the cops can move at most $\sigma / \psi \cdot 17\rho = \rho/2$.

\begin{figure}[h]
  \centering
  \includegraphics[width=0.8\textwidth]{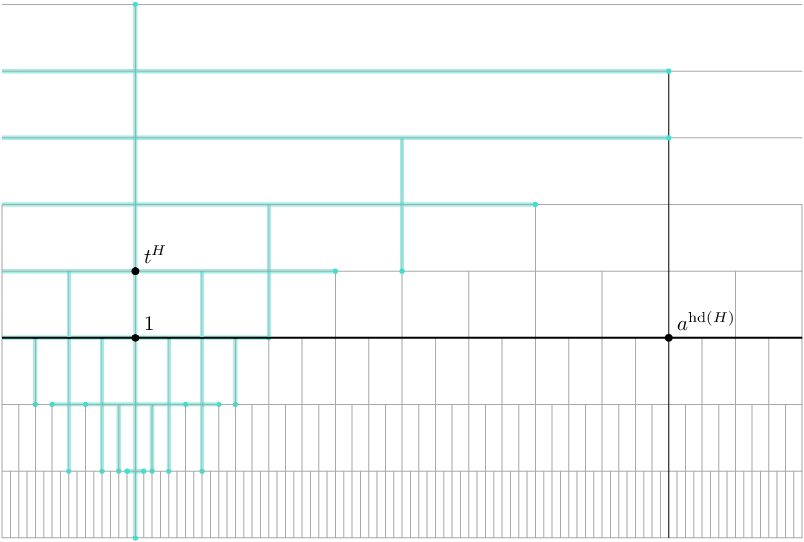} 
  \caption{A ball of radius $r=4$ centered at $t^H$ in the standard sheet for $H=1$. The largest $M$ such that the ball reaches $a^Mt^x$ for some $x\in \Z$ is $M=\operatorname{hd}(H) = 2^{
  r -1 +H} = 16$.}
  \label{fig:ball_radius_4}
\end{figure}

    We will show that the robber can choose one of the paths and follow along it without being caught by a cop.
    Let
    $$
    \operatorname{hd}(H) := \max \left\{ M \colon d(t^m, a^Mt^x) \leq 2\rho ,  m \leq H, x \in \Z \right\}
    $$
    be the maximal horizontal distance (as measured projected down on the axis) of a cop with height $m \leq H$ that can reach distance $2\rho$ (by first walking up to $\rho$ during the meta-step and then actually reaching $\rho$), see Figure \ref{fig:ball_radius_4}. Note that the best way to reach far horizontally is to first reach up by distance $2\rho-1$ and then over by distance $1$. This corresponds to a horizontal distance of $(n+1)^{2\rho-1}$ at the level where the cop is positioned and to a horizontal distance of at most $(n+1)^{3\rho-1}$ on the axis for a cop at height at most $\rho/2 + \sigma \leq \rho$. Therefore $ \operatorname{hd}(\rho) < (n+1)^{3\rho}$. Since the points $a^{k(n+1)^{8\rho} + i}$ with $0 \leq i \leq n$ lie more than $n+1$ times the distance $(n+1)^{3\rho}$ apart from each other (for different $k$), every cop can block at most one of the $n+1$ paths. There are only $n$ cops, so the robber can choose one. Thus the cops can never capture the robber and they can also not protect the ball, since the robber always stays within the rectangles, and in particular within the radius $R$ around the neutral element.
\end{proof}

\begin{remark}
    It is possible to define a quasi-retraction from $\operatorname{BS}(1,m)$ to a sheet, but since the sheet is hyperbolic and hence has strong cop number $1$, one can not use the quasi-retraction Theorem \ref{thm:retract} to infer Theorem \ref{thm:scop-bs}.
\end{remark}

\subsection{Grigorchuk's group}\label{sec:grigorchuk}

In this section we prove that the cop numbers of direct products of infinite groups are infinite. We apply this result to Grigorchuk's group and other branch groups. We first notice that the following generalization of Theorem \ref{thm:scop-Zn} holds.

\begin{theorem}\label{thm:scop-N2}
    The quadrant $\mathbb{N} \times \mathbb{N} \subseteq \mathbb{Z}^2$ has infinite strong cop number.
\end{theorem}
\begin{proof}
    After replacing $\mathbb{Z}^2$ by $\mathbb{N}^2$, the proofs of \cite[Theorem C]{LMR23} and Theorem \ref{thm:scop-Zn} can be followed word by word. 

		Alternatively, once one has defined the winning strategies for the robber in the weak cop game, the result follows from \Cref{thm:meta-gaming}.
\end{proof}

The idea for the following theorem was communicated to us by Francesco Fournier-Facio.

\begin{theorem}\label{thm:scop-products}
    If $G_1, \ldots, G_k$ with $k\geq 2$ are finitely generated infinite groups and $G=G_1 \times \ldots \times G_k$, then $\sCop(G) = \infty$.
\end{theorem}
\begin{proof}
A Cayley graph $\Gamma$ of an infinite finitely generated group has unbounded diameter. There exists a geodesic ray $\gamma \colon [0, \infty) \to \Gamma$ starting at the neutral element $\gamma(0) = e$. The word length $\ell \colon \Gamma \to \mathbb{Z}_{\geq 0}$ can be used to define a retraction\footnote{This is a retraction in the sense of metric spaces, but not in the sense of graph retractions as in \cite[Example 4.1(2)]{LMR23}. However, it is a quasi-retraction.} $\gamma \circ \ell \colon \Gamma \to \operatorname{im}(\gamma)$: indeed, for all $g,h \in \Gamma$ 
$$
d(\gamma(\ell(g)), \gamma(\ell(h))) = \lvert \ell(h) - \ell(g)\rvert = \vert d(e,h) - d(e,g) \rvert \leq d(g,h)
$$
by the triangle inequality. 
Let $\Gamma_1, \Gamma_2$ be the Cayley graphs of $G_1,G_2$ and let $\gamma_1, \gamma_2 $ and $\ell_1,\ell_2$ as above. 
In the Cayley graph of $G$, we obtain a retraction 
$r := \gamma_1\circ\ell_1 \times \gamma_2 \circ \ell_2\colon G \twoheadrightarrow G_1 \times G_2 \to \operatorname{im}(\gamma_1) \times \operatorname{im}(\gamma_2)$. The codomain is isometric to $\mathbb{Z}_{>0} \times \mathbb{Z}_{>0}$. By Theorems \ref{thm:retract} and \ref{thm:scop-N2}, $\sCop(G_1\times \ldots \times G_k) = \infty$.
\end{proof}

The Grigorchuk group and branch groups $G$ more generally contain a finite index subgroup $H$ that is isomorphic to a direct product $L\times \ldots \times L$ with at least two factors, where $L$ is a finitely generated infinite group (since $H$ is finite index in the finitely generated infinite group $G$) \cite[Definition 1.1]{BGZ03}.

\begin{theorem}\label{thm:scop-Grigorchuk}
    The Grigorchuk group and, more generally, finitely generated branch groups have infinite strong cop numbers.
\end{theorem}
\begin{proof}
    It is clear that branch groups $G$ are infinite. Any finite index subgroup $H$ that is isomorphic to a direct product $L\times  \ldots \times L$ with at least two finitely generated infinite factors has $\sCop(H) = \infty$ by Theorem \ref{thm:scop-products}. Finite index subgroups have quasi-isometric Cayley graphs, so $\sCop(G)=\infty$.
\end{proof}

\subsection{Thompson's group $F$}\label{sec:thompsonF}

It is well known among experts that Thompson's group $F$ admits a retraction to $\mathbb{Z}^2$, see for instance \cite[Cor. 3.10]{BFZ24arxiv} or \cite[Section 4]{CoMa24arxiv}. By Theorems \ref{thm:retract} and \ref{thm:scop-Zn}, then $\sCop(F) \geq \sCop(\mathbb{Z}^2) = \infty $. The analogous proof for the weak cop number has previously been used for the weak cop number by \cite{CoMa24arxiv}.

\begin{theorem}\label{thm:scopF=infty}
    The strong cop number of Thompson's group $F$ is infinite.
\end{theorem}

\section{Intermediate cop numbers}\label{sec:intermediate}

The main remaining open question is that of intermediate cop numbers. Partial results are known for graphs that are not Cayley-graphs of groups \cite[Corollary O]{LMR23}, but the question is wide open for finitely generated groups.

\begin{customquestionA}[\cite{LMR23} Question K]\label{q:intermediate}
Are there any finitely generated groups with weak and/or strong cop numbers not equal to $1$ or $\infty$?
\end{customquestionA}

A finitely generated group that acts geometrically on a CAT(0)-space is called a CAT(0)-group.

\begin{proposition}\label{prop:intermediate-CAT0}
    If $G$ is a CAT(0)-group, then $\sCop(G) \in \{ 1,\infty \} $.
\end{proposition}
\begin{proof} By Theorem \ref{thm:scop-hyperbolic} we may assume that $G$ is a finitely generated non-hyperbolic group $G$ that acts geometrically on a CAT(0)-space $X$. By \cite[Theorem 3.1, Chapter III.$\Gamma$]{BH99}, $X$ is known to contain an isometrically embedded copy of the Euclidean plane. Since CAT(0)-spaces admit retractions to convex subspaces, \Cref{thm:retract} can be applied, to obtain $\sCop(G)\geq \sCop(\mathbb{Z}^2) = \infty$ by Theorem \ref{thm:scop-Zn}.
\end{proof}

The cop numbers are closely related to phenomena of non-positive curvature. The question might therefore be especially interesting to solve for related classes of groups such as acylindrically hyperbolic groups, hierarchically hyperbolic groups, Helly groups or lacunary hyperbolic groups. As finitely presented groups are quotients of CAT(0)-groups, special attention should be given to groups that are not finitely presented. 

Our investigation into the Baumslag-Solitar groups, Lamplighter groups and Grigorchuk's group points towards the possibility that groups with intermediate cop numbers are rare. We expect the following two results to hold, but do not have an opinion on the answer to Question \ref{qA:intermediate}.

\begin{conjecture}
    Hyperbolic groups have weak cop number $1$ or $\infty$.\footnote{This has been shown to be true in \cite{Leh25arxiv} during the publication process of this paper.}
\end{conjecture}
\begin{conjecture}
    Finitely presented groups have cop numbers $1$ or $\infty$.
\end{conjecture}

\begin{figure}
  \centering
  \includegraphics[width=0.6\textwidth]{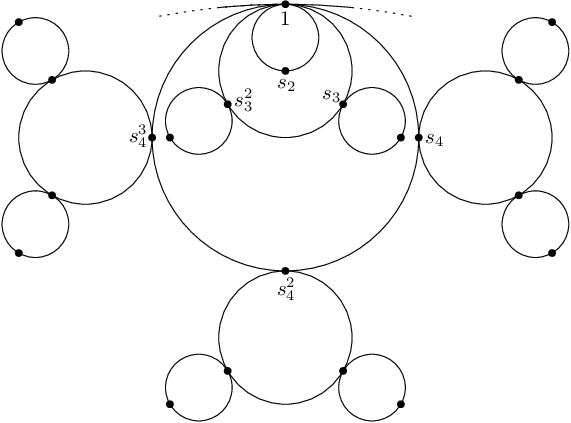} 
  \caption{The Cayley graph of the infinitely generated group $G = \langle s_2, s_3 \ldots \colon (s_n)^n = 1 \rangle$ with respect to the generating set $S = \{s_2, s_3, \ldots \}$ has weak and strong cop numbers $2$.}
  \label{fig:intermediate_infgen}
\end{figure}

\begin{remark}
Among infinitely generated groups, one can consider the Cayley graph of the group $G = \langle s_2, s_3, \ldots \colon (s_n)^n = 1 \rangle$ with respect to the generating set $S = \{s_2, s_3, \ldots \}$ to obtain a graph $\Gamma$ with $\sCop(\Gamma) = \wCop(\Gamma) = 2$, see Figure \ref{fig:intermediate_infgen}. However, for infinitely generated groups, the Cayley graphs for different generating sets may not be quasi-isometric, so the cop numbers may not be invariants of the group either.
\end{remark}

\newpage

\bibliographystyle{alpha}
\bibliography{references}

\end{document}